\documentclass[12pt,a4paper,svgnames]{amsart}
\usepackage[a4paper,left=2.25cm,right=2.25cm,top=3cm,bottom=3cm,headsep=1cm]{geometry}
\usepackage{mathdots}
\usepackage{tikz}\usetikzlibrary{decorations.markings,decorations.pathreplacing,arrows,matrix,positioning,shapes.multipart}
\colorlet{mygreen}{green!80!black}
\colorlet{myred}{red!90!black}
\definecolor{mygray}{rgb}{1.0,.96,.98}
\makeatletter

\newif\ifpgfcirclecrosssplitcustomfill
\tikzset{%
  circle cross split part fill/.code=%
    \def\pgf@lib@sh@ccs@list@fill{#1}\pgfcirclecrosssplitcustomfilltrue,%
  circle cross split uses custom fill/.is if=pgfcirclecrosssplitcustomfill}
\pgfdeclareshape{circle cross split}{%
  \nodeparts{text,two,three,four}%
  \savedanchor\centerpoint{%
    \pgfmathsetlength\pgf@xa{\pgfkeysvalueof{/pgf/inner xsep}}%
    \pgfmathsetlength\pgf@ya{\pgfkeysvalueof{/pgf/inner ysep}}%
    \pgf@x =\wd\pgfnodeparttextbox
    \pgf@yb=\dp\pgfnodeparttextbox
    \pgf@y=\dp\pgfnodeparttwobox
    \ifdim\pgf@yb>\pgf@y
      \pgf@y=\pgf@yb
    \fi
    \advance\pgf@x\pgf@xa
    \advance\pgf@y-\pgf@ya
    \advance\pgf@x.5\pgflinewidth
    \advance\pgf@y-.5\pgflinewidth
  }%
  \savedanchor\twoanchor{%
    \pgfmathsetlength\pgf@xa{\pgfkeysvalueof{/pgf/inner xsep}}%
    \pgfmathsetlength\pgf@ya{\pgfkeysvalueof{/pgf/inner ysep}}%
    \advance\pgf@x.5\pgflinewidth
    \advance\pgf@x\pgf@xa
    \advance\pgf@y.5\pgflinewidth
    \advance\pgf@y\pgf@ya
    \pgf@yb\dp\pgfnodeparttextbox
    \pgf@yc\dp\pgfnodeparttwobox
    \ifdim\pgf@yb>\pgf@yc
      \pgf@yc\pgf@yb
    \fi
    \advance\pgf@y\pgf@yc
  }%
  \savedanchor\threeanchor{%
    \pgfmathsetlength\pgf@ya{\pgfkeysvalueof{/pgf/inner ysep}}%
    \pgf@x\wd\pgfnodeparttextbox
    \pgf@yb\dp\pgfnodeparttextbox
    \pgf@yc\dp\pgfnodeparttwobox
    \ifdim\pgf@yb>\pgf@yc
      \pgf@yc\pgf@yb
    \fi
    \advance\pgf@y-\pgf@yc
    \advance\pgf@y-2\pgf@ya
    \advance\pgf@y-\pgflinewidth
    \pgf@yb\ht\pgfnodepartthreebox
    \pgf@yc\ht\pgfnodepartfourbox
    \ifdim\pgf@yb>\pgf@yc
      \pgf@yc\pgf@yb
    \fi
    \advance\pgf@y-\pgf@yc
    \advance\pgf@x-\wd\pgfnodepartthreebox
  }%
  \savedanchor\fouranchor{%
    \pgfmathsetlength\pgf@xa{\pgfkeysvalueof{/pgf/inner xsep}}%
    \advance\pgf@x\wd\pgfnodepartthreebox
    \advance\pgf@x2\pgf@xa
    \advance\pgf@x\pgflinewidth
  }%
  \saveddimen\radius{%
    \pgf@y\ht\pgfnodeparttextbox
    \pgf@yb\ht\pgfnodeparttwobox
    \ifdim\pgf@yb>\pgf@y
      \pgf@y\pgf@yb
    \fi
    \pgf@yc\dp\pgfnodeparttextbox
    \pgf@yb\dp\pgfnodeparttwobox
    \ifdim\pgf@yc>\pgf@yb
      \advance\pgf@y\pgf@yc
    \else
      \advance\pgf@y\pgf@yb
    \fi
    \pgf@yb\ht\pgfnodepartthreebox
    \ifdim\pgf@yb<\ht\pgfnodepartfourbox
      \pgf@yb\ht\pgfnodepartfourbox
    \fi
    \pgf@yc\dp\pgfnodepartthreebox
    \ifdim\pgf@yc<\dp\pgfnodepartfourbox
      \advance\pgf@yb\dp\pgfnodepartfourbox
    \else
      \advance\pgf@yb\pgf@yc
    \fi
    \ifdim\pgf@yc>\pgf@y
      \pgf@y\pgf@yc
    \fi
    \pgfmathsetlength\pgf@ya{\pgfkeysvalueof{/pgf/inner ysep}}%
    \advance\pgf@y2\pgf@ya
    \pgf@x\wd\pgfnodeparttextbox
    \pgf@xa\wd\pgfnodepartthreebox
    \pgf@xb\wd\pgfnodeparttwobox
    \pgf@xc\wd\pgfnodepartfourbox
    \ifdim\pgf@xa>\pgf@x
      \pgf@x\pgf@xa
    \fi
    \ifdim\pgf@xb>\pgf@x
      \pgf@x\pgf@xb
    \fi
    \ifdim\pgf@xc>\pgf@x
      \pgf@x\pgf@xc
    \fi
    \pgfmathsetlength\pgf@xa{\pgfkeysvalueof{/pgf/inner xsep}}%
    \advance\pgf@x2\pgf@xa
    \ifdim\pgf@y>\pgf@x
      \pgf@x\pgf@y
    \fi
    \advance\pgf@x.5\pgflinewidth
    \pgfmathsetlength{\pgf@xb}{\pgfkeysvalueof{/pgf/minimum width}}%
    \pgfmathsetlength{\pgf@yb}{\pgfkeysvalueof{/pgf/minimum height}}%
    \ifdim\pgf@x<.5\pgf@xb
        \pgf@x=.5\pgf@xb
    \fi
    \ifdim\pgf@x<.5\pgf@yb
        \pgf@x=.5\pgf@yb
    \fi
    \pgfmathsetlength{\pgf@xb}{\pgfkeysvalueof{/pgf/outer xsep}}%
    \pgfmathsetlength{\pgf@yb}{\pgfkeysvalueof{/pgf/outer ysep}}%
    \ifdim\pgf@xb<\pgf@yb
      \advance\pgf@x\pgf@yb
    \else
      \advance\pgf@x\pgf@xb
    \fi
  }%
  \inheritanchorborder[from=circle]%
  \inheritanchor[from=circle]{north}%
  \inheritanchor[from=circle]{north west}%
  \inheritanchor[from=circle]{north east}%
  \inheritanchor[from=circle]{center}%
  \inheritanchor[from=circle]{west}%
  \inheritanchor[from=circle]{east}%
  \inheritanchor[from=circle]{mid}%
  \inheritanchor[from=circle]{mid west}%
  \inheritanchor[from=circle]{mid east}%
  \inheritanchor[from=circle]{base}%
  \inheritanchor[from=circle]{base west}%
  \inheritanchor[from=circle]{base east}%
  \inheritanchor[from=circle]{south}%
  \inheritanchor[from=circle]{south west}%
  \inheritanchor[from=circle]{south east}%
  \anchor{two}{\twoanchor}%
  \anchor{three}{\threeanchor}%
  \anchor{four}{\fouranchor}%
  \inheritbackgroundpath[from=circle]
  \beforebackgroundpath{%
    \pgfutil@tempdima=\radius
    \pgfmathsetlength{\pgf@xb}{\pgfkeysvalueof{/pgf/outer xsep}}%
    \pgfmathsetlength{\pgf@yb}{\pgfkeysvalueof{/pgf/outer ysep}}%
    \ifdim\pgf@xb<\pgf@yb
      \advance\pgfutil@tempdima by-\pgf@yb
    \else
      \advance\pgfutil@tempdima by-\pgf@xb
    \fi
    \advance\pgfutil@tempdima by-.5\pgflinewidth%
    \pgfsetshortenstart{0pt}%
    \pgfsetshortenend{0pt}%
    \pgfsetarrows{-}%
    \pgfpathmoveto
      {\pgfpointadd{\centerpoint}{\pgfqpoint{-\pgfutil@tempdima}{0pt}}}%
    \pgfpathlineto
      {\pgfpointadd{\centerpoint}{\pgfqpoint{\pgfutil@tempdima}{0pt}}}%
    \pgfpathmoveto
      {\pgfpointadd{\centerpoint}{\pgfqpoint{0pt}{-\pgfutil@tempdima}}}%
    \pgfpathlineto
      {\pgfpointadd{\centerpoint}{\pgfqpoint{0pt}{\pgfutil@tempdima}}}%
    \pgfusepathqstroke
  }%
  \behindbackgroundpath{%
    \pgfutil@tempdima=\radius
    \pgfmathsetlength{\pgf@xb}{\pgfkeysvalueof{/pgf/outer xsep}}%
    \pgfmathsetlength{\pgf@yb}{\pgfkeysvalueof{/pgf/outer ysep}}%
    \ifdim\pgf@xb<\pgf@yb
      \advance\pgfutil@tempdima by-\pgf@yb
    \else
      \advance\pgfutil@tempdima by-\pgf@xb
    \fi
    \advance\pgfutil@tempdima by-.5\pgflinewidth%
    \ifpgfcirclecrosssplitcustomfill%
      \pgf@lib@sh@rs@process@list{\pgf@lib@sh@ccs@list@fill}{4}%
      {%
        \pgfmathloop
           \ifnum\pgfmathcounter>4%
           \else%
             \pgf@lib@sh@getalpha\pgf@lib@sh@rs@number{\pgfmathcounter}%
              \edef\pgf@tempa
                {\csname pgf@lib@sh@rs@\pgf@lib@sh@rs@number @item\endcsname}%
              \ifx\pgf@tempa\pgf@lib@sh@rs@nonetext\else
                \pgfsetfillcolor{\pgf@tempa}%
                \pgf@lib@sh@ccs@angles{\pgfmathcounter}%
                \pgfpathmoveto{\centerpoint}%
                \pgfpathlineto{\pgfpointadd{\centerpoint}
                  {\pgfqpointpolar{\pgf@lib@sh@ccs@angle}{\pgfutil@tempdima}}}%
                \pgfpatharc{\pgf@lib@sh@ccs@angle}{\pgf@lib@sh@ccs@angle@}
                  {\pgfutil@tempdima}%
                \pgfpathclose
                \pgfusepathqfill
              \fi
        \repeatpgfmathloop
      }%
    \fi}}
\def\pgf@lib@sh@ccs@angles#1{%
  \ifcase#1\or\def\pgf@lib@sh@ccs@angle{90}%
           \or\def\pgf@lib@sh@ccs@angle{0}%
           \or\def\pgf@lib@sh@ccs@angle{180}%
           \else\def\pgf@lib@sh@ccs@angle{270}%
  \fi
  \edef\pgf@lib@sh@ccs@angle@{\number\numexpr\pgf@lib@sh@ccs@angle+90\relax}}
\makeatother
\tikzset{math mode/.style = {execute at begin node=$, execute at end node=$}}
\tikzset{sixvin/.style={>=triangle 45,decoration={markings,mark = at position #1 with { \arrow{>} }},postaction={decorate}},sixvin/.default=0.6}
\tikzset{sixvout/.style={>=triangle 45,decoration={markings,mark = at position #1 with { \arrow{<} }},postaction={decorate}},sixvout/.default=0.6}
\tikzset{arrow/.style={postaction={decorate,decoration={markings,mark = at position #1 with {\arrow{Straight Barb[line width=0.2mm,length=1.2mm]}}}}},arrow/.default=0.5}
\tikzset{invarrow/.style={postaction={decorate,decoration={markings,mark = at position #1 with {\arrowreversed{Straight Barb[line width=0.2mm,length=1.2mm]}}}}},invarrow/.default=0.5}
\usepackage{hyperref}
\newtheorem{theorem}{Theorem}
\allowdisplaybreaks[1]

\renewcommand\ss{\scriptstyle}

\newcommand\ZZ{{\mathbb Z}}
\newcommand\rem[2][]{}
\title{Integrability and combinatorics}
\author{Paul Zinn-Justin}
\address{Paul Zinn-Justin, School of Mathematics and Statistics, The University of Melbourne, 
Victoria 3010, Australia}
\email{pzinn@unimelb.edu.au}
\thanks{The author thanks A.~Gunna for a careful reading of the manuscript.}
\begin{document}
\begin{abstract}
We discuss the use of methods coming from integrable systems to study problems of enumerative and algebraic combinatorics, and develop two examples:
the enumeration of Alternating Sign Matrices and related combinatorial objects, and the theory of symmetric polynomials.
\end{abstract}

\maketitle



\section{Introduction}
\subsection{Generalities}
Combinatorics and mathematical physics have many points of contact. Here, we focus on a specific form of interaction, which is the use
of methods coming from quantum integrable systems to solve problems of a combinatorial nature. A typical application is
to enumerative combinatorics. As the field of combinatorics is expanding rapidly,
problems of enumeration become more and more difficult and direct combinatorial proofs are often extremely
complicated and tedious. This is where physical ideas, and in particular integrability,
can come in to provide conceptually simpler proofs. In its most basic form, one may hope that integrability,
under the guise of exactly solvability of lattice of models of statistical mechanics,
allows us to perform exact computations of partition
functions and therefore enumerate the underlying combinatorial objects.
Going beyond this naive viewpoint, one observes that
there is a deeper connection between integrability and combinatorics, and more specifically
algebraic combinatorics. In this short review, we can only give hints of this connection, and of the shared underlying algebraic structures
and representation theory.

\subsection{Plan}
In what follows, we shall give two examples, one for each of the two connections outlined above, where
these ideas turned out particularly fruitful:
\begin{itemize}
\item In \S\ref{sec:6v}, we discuss the connection between the Six-Vertex model with Domain Wall Boundary Conditions,
and various enumerative problems, in particular Alternating Sign Matrices.
\item In \S\ref{sec:sym}, we reinterpret the theory of symmetric polynomials, a classical topic of algebraic combinatorics,
  in terms of quantum integrable systems, and show some applications, focusing on the prototypical
  case of Schur polynomials.
\end{itemize}
\rem{add ref to other chapter about QIS. }

\subsection{The Six-Vertex model}
In both these examples, a key role is played by
the {\em six vertex model}, an important model of classical statistical
mechanics in two dimensions. It first appeared as a model for (two-dimensional) ice, which
was solved by Lieb \cite{Lieb1} in 1967 by means of Bethe Ansatz, followed by several 
generalizations \cite{Lieb2, Lieb3, Lieb4,Suth}.

\begin{figure}
\begin{tikzpicture}
\matrix[cells={scale=.8},column sep=1cm,row sep=0.2cm] {
\draw[sixvin] (-1,0) -- (0,0); \draw[sixvin] (0,-1) -- (0,0); \draw[sixvout] (1,0) -- (0,0); \draw[sixvout] (0,1) -- (0,0);
&
\draw[sixvout] (-1,0) -- (0,0); \draw[sixvin] (0,-1) -- (0,0); \draw[sixvin] (1,0) -- (0,0); \draw[sixvout] (0,1) -- (0,0);
&
\draw[sixvin] (-1,0) -- (0,0); \draw[sixvout] (0,-1) -- (0,0); \draw[sixvin] (1,0) -- (0,0); \draw[sixvout] (0,1) -- (0,0);
\\
\node{$a_1$};&\node{$b_1$};&\node{$c_1$};
\\[2mm]
\draw[sixvout] (-1,0) -- (0,0); \draw[sixvout] (0,-1) -- (0,0); \draw[sixvin] (1,0) -- (0,0); \draw[sixvin] (0,1) -- (0,0);
&
\draw[sixvin] (-1,0) -- (0,0); \draw[sixvout] (0,-1) -- (0,0); \draw[sixvout] (1,0) -- (0,0); \draw[sixvin] (0,1) -- (0,0);
&
\draw[sixvout] (-1,0) -- (0,0); \draw[sixvin] (0,-1) -- (0,0); \draw[sixvout] (1,0) -- (0,0); \draw[sixvin] (0,1) -- (0,0);
\\
\node{$a_2$};&\node{$b_2$};&\node{$c_2$};
\\
};
\end{tikzpicture}
\caption{Local configurations of the six-vertex model.}
\label{fig:6vweights}
\end{figure}
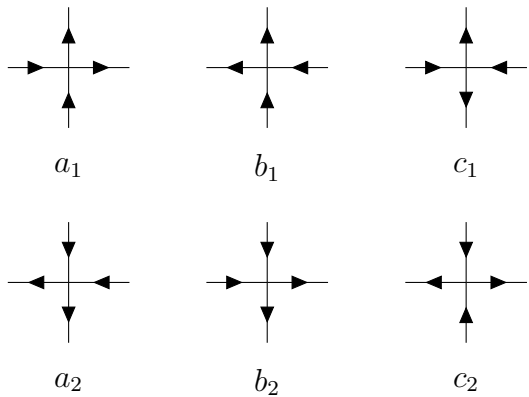

The six-vertex model is a statistical model defined on a (subset of the) 
square lattice; its configurations are obtained by putting arrows (two possible directions) on each edge of the lattice,
with the additional rule
that at each vertex, there are as many incoming arrows as outgoing ones.
Around a given vertex, there are only 6 configurations of edges which respect this ``arrow conservation'' rule, see Fig.~\ref{fig:6vweights}, hence the name of the model.

To each such a configuration one associates a Boltzmann weight which is a product of local weights for each vertex,
leading to six parameters $a_1,a_2,b_1,b_2,c_1,c_2$ as on Fig.~\ref{fig:6vweights}.
The partition function is then defined as
\[
Z=\sum_{\text{configurations}}\,
\prod_{x\in\{a_1,a_2,b_1,b_2,c_1,c_2\}}
x^{\#\text{vertices of type }x}
\]

\section{The Six-Vertex model with Domain Wall Boundary Conditions}\label{sec:6v}
Domain Wall Boundary Conditions (DWBC) are specific boundary conditions for the Six-Vertex model \cite{Kor}
which turned out to possess remarkable properties
\cite{Iz-6V, ICK}, as will be discussed below.

The connection to combinatorics appeared in 1996 in the work of Kuperberg \cite{Kup-ASM} who noticed that configurations
of the Six-Vertex model with DWBC are in bijection with Alternating Sign Matrices,
a famous object in combinatorics \cite{Bressoud}, and used this bijection
to give a simple proof of the Alternating Sign Matrix
conjecture, that is to solve the underlying enumeration problem. We shall give a modified version of his result below.

In the 2000s, the Six-Vertex model with Domain Wall Boundary Conditions became popular again as a simple model
to test sensitivity of lattice models to boundary conditions \cite{artic12}. In particular, it exhibits
in the limit of large size {\em limiting shapes}, i.e., spatial phase separation induced by the boundary, 
see \cite{Stephan-tilings} and references therein.

\subsection{Definition}
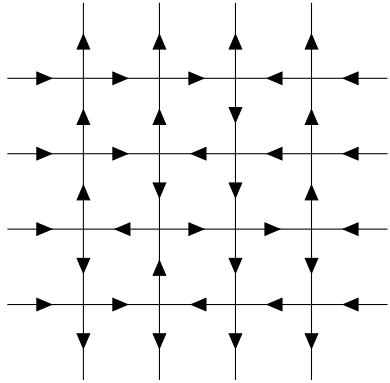
\begin{figure}
\begin{tikzpicture}
\draw[sixvin] (0,1) -- (1,1);\draw[sixvin] (1,1) -- (2,1);\draw[sixvin] (5,1) -- (4,1);\draw[sixvin] (4,1) -- (3,1);\draw[sixvin] (3,1) -- (2,1);\draw[sixvin] (0,2) -- (1,2);\draw[sixvin] (2,2) -- (1,2);\draw[sixvin] (2,2) -- (3,2);\draw[sixvin] (3,2) -- (4,2);\draw[sixvin] (5,2) -- (4,2);\draw[sixvin] (0,3) -- (1,3);\draw[sixvin] (1,3) -- (2,3);\draw[sixvin] (3,3) -- (2,3);\draw[sixvin] (4,3) -- (3,3);\draw[sixvin] (5,3) -- (4,3);\draw[sixvin] (0,4) -- (1,4);\draw[sixvin] (1,4) -- (2,4);\draw[sixvin] (2,4) -- (3,4);\draw[sixvin] (4,4) -- (3,4);\draw[sixvin] (5,4) -- (4,4);
\draw[sixvin] (1,1) -- (1,0);\draw[sixvin] (1,2) -- (1,1);\draw[sixvin] (1,2) -- (1,3);\draw[sixvin] (1,3) -- (1,4);\draw[sixvin] (1,4) -- (1,5);\draw[sixvin] (2,1) -- (2,0);\draw[sixvin] (2,1) -- (2,2);\draw[sixvin] (2,3) -- (2,2);\draw[sixvin] (2,3) -- (2,4);\draw[sixvin] (2,4) -- (2,5);\draw[sixvin] (3,1) -- (3,0);\draw[sixvin] (3,2) -- (3,1);\draw[sixvin] (3,3) -- (3,2);\draw[sixvin] (3,4) -- (3,3);\draw[sixvin] (3,4) -- (3,5);\draw[sixvin] (4,1) -- (4,0);\draw[sixvin] (4,2) -- (4,1);\draw[sixvin] (4,2) -- (4,3);\draw[sixvin] (4,3) -- (4,4);\draw[sixvin] (4,4) -- (4,5);
\end{tikzpicture}
\caption{An example of configuration of the six-vertex model with Domain Wall Boundary Conditions.}
\label{fig:6vdwbc}
\end{figure}

Domain Wall Boundary Conditions (DWBC) consist in considering
the six-vertex model on a
$n\times n$ square domain, and fixing all external edges: vertical (resp.\ horizontal)
external edges are fixed to be outgoing (resp.\ incoming). See Fig.~\ref{fig:6vdwbc} for a $n=4$ example.
We denote $\text{DWBC}_n$ the set of such configurations, and $Z_n$ the corresponding partition
function.

It is not hard to show (and the reader is encouraged to check using
any of the alternative representations of $\text{DWBC}_n$ given below)
that the numbers of vertices of type $a_1$ and $a_2$ are equal in a DWBC configuration,
and similarly for $b_1$ and $b_2$; and that there are exactly $n$ vertices of type $c_1$ more than there are of 
type $c_2$ (one extra vertex $c_1$ per row or per column). Therefore, the partition function can also be written
\[
Z_n=
(c_1/c_2)^{n/2}
\sum_{\mathcal C\in\text{DWBC}_n}\,
a^{\#\text{vertices of type }a}
b^{\#\text{vertices of type }b}
c^{\#\text{vertices of type }c}
\]
where $a=\sqrt{a_1a_2}$, $b=\sqrt{b_1b_2}$, $c=\sqrt{c_1c_2}$,
and a vertex of type $a$ is a vertex of type either $a_1$ or $a_2$, and similarly for $b$ and $c$.

\subsection{Alternate representations}\label{ssec:alt}
We pause to discuss several mappings between configurations of the six-vertex model with DWBC
and other interesting combinatorial objects, as well as their interrelations.

\subsubsection{Lattice paths}\label{ssec:lp}
Let us begin with a trivial bijection: let us relabel edges of the lattice such that edges with right or up
arrows are ``occupied'', whereas those with left or down arrows are ``empty''. The configuration of
Figure~\ref{fig:6vdwbc} is redrawn this way on Figure~\ref{fig:lp}. This shows that the six-vertex model
can also be viewed as a model of {\em lattice paths} going North/East, which can touch but not
cross each other (such paths are also known as ``osculating walkers'') \cite{Brak-ASM}.
The DWBC simply mean that the endpoints
of the paths lie on the West and North sides of the square lattice.

\tikzset{part/.style={circle,thin,solid,draw=black,inner sep=1.2pt}}
\tikzset{occupied/.style={draw=none},empty/.style={thin,densely dotted}}
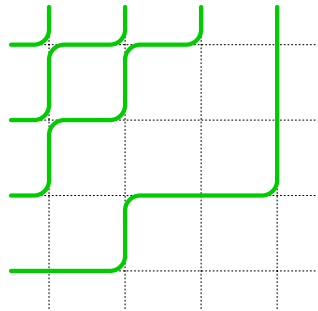
\begin{figure}
\begin{tikzpicture}[line cap=round]
\draw[empty] (0.5,0.5) grid (4.5,4.5);
\begin{scope}[ultra thick,mygreen,rounded corners=2mm]
\draw (0.5,1)  -- (2,1) -- (2,2) -- (4,2) -- (4,4.5) ;
\draw (0.5,2) -- (1,2) -- (1,3) -- (2,3) -- (2,4) -- (3,4) -- (3,4.5) ;
\draw (0.5,3) -- (1,3) -- (1,4) -- (2,4) -- (2,4.5) ;
\draw (0.5,4) -- (1,4) -- (1,4.5) ;
\end{scope}
\end{tikzpicture}
\caption{An example of a lattice path configuration.}
\label{fig:lp}
\end{figure}

\subsubsection{Rook placements}\label{ssec:rook}
In each row (or column) of the $n\times n$ grid, given a configuration in $\text{DWBC}_n$,
there must be at least one configuration of type $c$, because the boundary arrows are opposite.
There is therefore a natural subset of $\text{DWBC}_n$, namely
the configurations where there is a single type $c$ vertex per row (which implies the same property for columns);
the $c$ vertices are all of type $c_1$.
These configurations are in obvious bijection with (complete) ``rook placements'', i.e., 
non-attacking configurations of $n$ rooks
on a $n\times n$ chessboard (where the rooks sit at $c$ vertices). Equivalently, these configurations are in bijection
with permutations of $\mathcal S_n$, where recording row/columns of $c$ vertices produces a permutation. In particular,
there are $n!$ of them. We shall come back to these configurations shortly.

\subsubsection{Height functions}
A first nontrivial bijection, whose existence is directly related to the arrow conservation law that is built-in the definition of the
model, is to associate a height to each face of the square lattice, in such a way that when one goes from one face
to an adjacent one, the height varies by $+1$ (resp.\ $-1$) if the arrow on the edge separating them points
left (resp.\ right). This only determines the height up to an overall constant, which we fix by imposing that the
height at the top left of the lattice is zero. The same example of Fig.~\ref{fig:6vdwbc} is depicted as a height function
on Fig.~\ref{fig:height}.

\begin{figure}
\begin{tikzpicture}[math mode]
\matrix[matrix of nodes,column sep={1cm,between origins},row sep={1cm,between origins}]{
0&1&2&3&4
\\
1&2&3&2&3
\\
2&3&2&1&2
\\
3&2&3&2&1
\\
4&3&2&1&0
\\
};
\end{tikzpicture}
\caption{An example of height function.}
\label{fig:height}
\end{figure}

It is not hard to see that the set of height functions obtained this way from a DWBC configuration is exactly
\[
\text{H}_n=\left\{
(h_{ij})_{i,j=0,\ldots,n}\ \Bigg|\  
\begin{aligned}
&h_{i0}=i,\ h_{0j}=j,\ h_{in}=n-i,\ h_{nj}=n-j,&& i,j=0,\ldots,n\\
&h_{i+1j}-h_{ij},\ h_{ij+1}-h_{ij}\in \{-1,1\},&&i,j=0,\ldots,n-1
\end{aligned}
\right\}
\]
(where height functions are indexed as matrices), and that the mapping is bijective.

There is an additional structure on the set $\text{H}_n$ above (that was perhaps not so apparent in the original
six-vertex formulation): it is a {\em lattice}\/ (in the sense of ordered sets). That is, $\text{H}_n$ is a poset
-- the order relation is pointwise $\le$ --
such that any pair of configurations possesses an infimum and a supremum (pointwise maximum, minimum).

The subset of ``rook placements'', under this bijection, becomes a sub-poset of the poset of height functions, which
has a natural interpretation: it is the type $A$ {\em Bruhat poset}, that is, the symmetric group $\mathcal S_n$
endowed with its
Bruhat order. From this point of view, one can view the whole of $\text{H}_n$ as the {\em MacNeille completion}\/ of
the type $A$ Bruhat poset \cite{LS-Cox}.

\subsubsection{Alternating Sign Matrices}\label{ssec:ASM}
Permutations can be represented as permutation matrices (matrices with a single $1$ per row/column, and zero elsewhere);
in view of the above, it is natural to ask how one can reconstruct the permutation matrix from a configuration in
$\text{DWBC}_n$ of the rook placement type. This is easily done using height functions: indeed, defining
\begin{equation}\label{eq:hgt2asm}
w_{ij}=\frac{1}{2}(h_{i\,j-1}+h_{i-1\,j}-h_{ij}-h_{i-1\,j-1})\qquad i,j=1,\ldots,n
\end{equation}
leads to the desired permutation matrix.

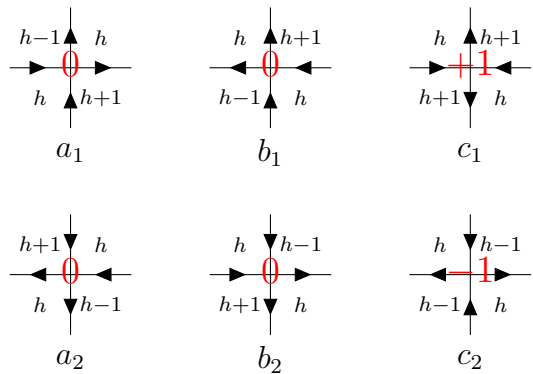
\begin{figure}
\begin{tikzpicture}
\matrix[asm/.style={red,yshift=.5mm,font=\large},cells={scale=.8},column sep=1cm,row sep=0.5cm] {
\draw[sixvin] (-1,0) -- (0,0); \draw[sixvin] (0,-1) -- (0,0); \draw[sixvout] (1,0) -- (0,0); \draw[sixvout] (0,1) -- (0,0);
\node[asm] {$0$};
\node at (-.5,.5) {$\ss h-1$};
\node at (.5,.5) {$\ss h$};
\node at (-.5,-.5) {$\ss h$};
\node at (.5,-.5) {$\ss h+1$};
&
\draw[sixvout] (-1,0) -- (0,0); \draw[sixvin] (0,-1) -- (0,0); \draw[sixvin] (1,0) -- (0,0); \draw[sixvout] (0,1) -- (0,0);
\node[asm] {$0$};
\node at (-.5,.5) {$\ss h$};
\node at (.5,.5) {$\ss h+1$};
\node at (-.5,-.5) {$\ss h-1$};
\node at (.5,-.5) {$\ss h$};
&
\draw[sixvin] (-1,0) -- (0,0); \draw[sixvout] (0,-1) -- (0,0); \draw[sixvin] (1,0) -- (0,0); \draw[sixvout] (0,1) -- (0,0);
\node[asm] {$+1$};
\node at (-.5,.5) {$\ss h$};
\node at (.5,.5) {$\ss h+1$};
\node at (-.5,-.5) {$\ss h+1$};
\node at (.5,-.5) {$\ss h$};
\\[-.5cm]
\node{$a_1$};&\node{$b_1$};&\node{$c_1$};
\\
\draw[sixvout] (-1,0) -- (0,0); \draw[sixvout] (0,-1) -- (0,0); \draw[sixvin] (1,0) -- (0,0); \draw[sixvin] (0,1) -- (0,0);
\node[asm] {$0$};
\node at (-.5,.5) {$\ss h+1$};
\node at (.5,.5) {$\ss h$};
\node at (-.5,-.5) {$\ss h$};
\node at (.5,-.5) {$\ss h-1$};
&
\draw[sixvin] (-1,0) -- (0,0); \draw[sixvout] (0,-1) -- (0,0); \draw[sixvout] (1,0) -- (0,0); \draw[sixvin] (0,1) -- (0,0);
\node[asm] {$0$};
\node at (-.5,.5) {$\ss h$};
\node at (.5,.5) {$\ss h-1$};
\node at (-.5,-.5) {$\ss h+1$};
\node at (.5,-.5) {$\ss h$};
&
\draw[sixvout] (-1,0) -- (0,0); \draw[sixvin] (0,-1) -- (0,0); \draw[sixvout] (1,0) -- (0,0); \draw[sixvin] (0,1) -- (0,0);
\node[asm] {$-1$};
\node at (-.5,.5) {$\ss h$};
\node at (.5,.5) {$\ss h-1$};
\node at (-.5,-.5) {$\ss h-1$};
\node at (.5,-.5) {$\ss h$};
\\[-.5cm]
\node{$a_2$};&\node{$b_2$};&\node{$c_2$};
\\
};
\end{tikzpicture}
\caption{Local configurations of the six-vertex model, their height function and ASM formulation (the latter in red).}
\label{fig:6vweights2}
\end{figure}

Let us now apply this mapping to the whole of $\text{H}_n$; for convenience, we reproduce on Fig.~\ref{fig:6vweights2}
the mappings between the local configurations. Also see Fig.~\ref{fig:exasm} for our running example.

\begin{figure}
\begin{tikzpicture}[every node/.style={red,font=\large},math mode]
\node at (1,4) {0}; \node at (2,4) {0}; \node at (3,4) {1}; \node at (4,4) {0}; \node at (1,3) {0}; \node at (2,3) {1}; \node at (3,3) {0}; \node at (4,3) {0}; \node at (1,2) {1}; \node at (2,2) {-1};\node at (3,2) {0}; \node at (4,2) {1}; \node at (1,1) {0}; \node at (2,1) {1}; \node at (3,1) {0}; \node at (4,1) {0}; 
\end{tikzpicture}
\caption{An example of an ASM.}
\label{fig:exasm}
\end{figure}

We observe that besides the values $1$ for vertices of type $c_1$ and $0$ for vertices of type $a$ and $b$, another
value is possible: $-1$ for vertices of type $c_2$. In fact, it is easy to show that the mapping \eqref{eq:hgt2asm}
leads to a bijection of $\text{H}_n$ (and therefore, $\text{DWBC}_n$) with the set of {\em Alternating Sign Matrices}
\begin{equation}
\text{ASM}_n=\left\{
(w_{ij})_{i,j=1,\ldots,n}\ \Bigg|\  
\begin{aligned}
&w_{ij}\in \{0,\pm 1\}
\\
&\text{+1s and -1s alternate on each row and column,}
\\[-2mm]
&\text{starting and ending with $1$s}
\end{aligned}
\right\}
\end{equation}
Alternating Sign Matrices (ASMs) were introduced by Robbins and Rumsey \cite{RR-ASM} in the context of Dodgson's condensation method
for computing determinants. Following this method naturally leads to an expansion of the determinant as a sum over ASMs, except the coefficient of
non-permutation matrices turns out to be zero. This can be remedied by introducing a deformation parameter into the formula, leading to the notion
of {\em $\lambda$-determinant}, defined by:
\newcommand\ldet{\det{}\!_\lambda}
\begin{equation}\label{eq:ldet}
\ldet M=\sum_{A\in\text{ASM}(n)}
\lambda^{\nu(A)} (1+\lambda)^{\mu(A)}
\prod_{i,j=1}^n M_{ij}^{A_{ij}}
\end{equation}
where $\mu(A)$ is the number of $-1$s in $A$, and $\nu(A)=\sum_{\substack{1\le i\le i'\le n\\ 1\le j'<j\le n}}
A_{ij}A_{i'j'}$ is a generalization of the inversion number of a permutation. At $\lambda=-1$,
one recovers the usual determinant.

We note that via the bijection to $\text{DWBC}_n$,
$\mu(A)$ is nothing but the number of vertices of type $c_1$.
Similarly, one can show \cite{artic55}
that $\nu(a)$ is the number of vertices of type $a_1$, or of type $a_2$, that is to say half the total
number of vertices of type $a$.

Let us for example choose the constant matrix $M=I_n$ with $(I_n)_{ij}=1$ in \eqref{eq:ldet}:
we immediately conclude from the above that its $\lambda$-determinant coincides with the
DWBC partition function for a particular choice of Boltzmann weights, namely,
\begin{equation}\label{eq:ldet2}
\ldet I_n = Z_n(a=\sqrt{\lambda},b=1,c_1=1,c_2=1+\lambda)
\end{equation}
Note that this is not the most general choice
of six-vertex model parameters, even taking into account the freedom to rescale all weights by a constant.
In particular, a natural quantity to compute is the number of ASMs, that is, $\#\text{ASM}_n=\#\text{DWBC}_n$;
it is {\em not}\/ a special case of \eqref{eq:ldet2}. The determination of $\#\text{ASM}_n$ has a long and rich history,
recounted in \cite{Bressoud,BP-ASM}; here we merely point out that the number of ASMs is equal to the number
of (at least)
two other, {\em a priori}\/ unrelated, families of combinatorial objects: Descending Plane Partitions \cite{And-MacDPP} and
Totally Symmetric Self-Complementary Plane Partitions \cite{And-TSSCPP}. These equalities, and various refinements thereof,
can be proven using tools from quantum integrability, see e.g.~\cite{Zeil-refASM,CP-orthog2,artic45,artic55}.

\subsubsection{Monotone triangles}
A DWBC configuration is entirely determined by the state of all vertical (or horizontal) edges.
Let us therefore record for each row the subset of up-pointing arrows: we obtain this way a triangular array
known as {\em monotone triangle}\/ \cite{MRR-ASM} (or gog triangle \cite{Zeil-ASM}). On our running example, one finds
\[
\begin{tikzpicture}[yscale=0.7]
\node at (1,4.5) {1};
\node at (2,4.5) {2};
\node at (3,4.5) {3};
\node at (4,4.5) {4};
\node at (1.5,3.5) {1};
\node at (2.5,3.5) {2};
\node at (3.5,3.5) {4};
\node at (2,2.5) {1};
\node at (3,2.5) {4};
\node at (2.5,1.5) {2};
\end{tikzpicture}
\]

In general, monotone triangles of size $n$ are triangular arrays of integers of the form
\[
\begin{tikzpicture}
\matrix[matrix of math nodes,row sep={1cm,between origins},column sep={1cm,between origins}] at (0,-.5) {
1 &<& 2 &<& \cdots &<&  n-1 &<& n
\\
& a_{n-1,1}&<&a_{n-1,2}&<&\cdots&<&a_{n-1,n-1}
\\
&&\ddots&&&&\iddots
\\
&&&a_{2,1}&<&a_{2,2}
\\
&&&&a_{1,1}
\\
};
\node[rotate=45] at (3.5,1) {$\le$};
\node[rotate=45] at (1.5,1) {$\le$};
\node[rotate=45] at (-0.5,1) {$\le$};
\node[rotate=45] at (-2.5,1) {$\le$};
\node[rotate=-45] at (2.5,1) {$\le$};
\node[rotate=-45] at (0.5,1) {$\le$};
\node[rotate=-45] at (-1.5,1) {$\le$};
\node[rotate=-45] at (-3.5,1) {$\le$};
\node[rotate=45] at (2.5,0) {$\le$};
\node[rotate=-45] at (-2.5,0) {$\le$};
\node[rotate=45] at (1.5,-1) {$\le$};
\node[rotate=-45] at (-1.5,-1) {$\le$};
\node[rotate=45] at (0.5,-2) {$\le$};
\node[rotate=-45] at (-0.5,-2) {$\le$};
\end{tikzpicture}
\]
One can check once again that their set is in bijection with $\text{DWBC}_n$; it also has an obvious order relation (pointwise $\le$)
which is the same one that we have defined on $\text{H}_n$ up to the various bijections.

Monotone triangles are also called {\em strict Gelfand--Tsetlin patterns} (with fixed first row $12\ldots n$),
because if one replaces strict inequality
along rows with weak inequality, then one recovers the definition of Gelfand--Tsetlin patterns.

It is well-known that Gelfand--Tsetlin patterns are in bijection with {\em Semi-Standard Young Tableaux} (SSYT) (see \S\ref{ssec:schur}),
where the first row determines the shape of the tableau. Drawing tableaux with the French notation,
we obtain on our running example:
\newdimen{\cellsize}
\newcommand\bigboxes{\setlength{\cellsize}{18pt}\def\boxformat{}}
\newcommand\medboxes{\setlength{\cellsize}{14.22pt}\def\boxformat{}}
\newcommand\smallboxes{\setlength{\cellsize}{10pt}\def\boxformat{\scriptstyle}}
\newcommand\tinyboxes{\setlength{\cellsize}{6pt}\def\boxformat{\scriptscriptstyle}}
\medboxes
\tikzset{tableaubox/.style={draw=black,thin,sharp corners,solid,minimum size=\cellsize,inner sep=0pt}}
\tikzset{tableau/.style={matrix,name=tab,matrix anchor=tab-1-1.south west,inner sep=1pt,matrix of math nodes,cells={anchor=center,draw=black,thin,solid,arrows=-},nodes={tableaubox,execute at begin node=\boxformat},nodes in empty cells,row sep={\cellsize,between origins},column sep={\cellsize,between origins}}}
%
\newcommand\colorcell[1]{|[fill=#1]|}
\newcommand\graycell{|[fill=gray]|}
\newcommand\thickcell{|[line width=2pt,minimum size=\cellsize-1.6pt]|}
\newcommand\missingcell{|[draw=none]|}
%
\makeatletter
\newcommand\cellextra[1]{#1\expandafter\tikz@lib@matrix@start@cell}
\makeatother
\newcommand\hdotscell{\cellextra{\draw[dotted] (-0.5*\cellsize,0.5*\cellsize) --++(\cellsize,0);}\missingcell}
\newcommand\vdotscell{\cellextra{\draw[dotted] (-0.5*\cellsize,-0.5*\cellsize) --++(0,\cellsize);}\missingcell}
\newcommand\vhdotscell{\cellextra{\draw[dotted] (-0.5*\cellsize,-0.5*\cellsize) --++(0,\cellsize) --++(\cellsize,0);}\missingcell}
\newcommand\hcell{\cellextra{\draw (-0.5*\cellsize,0.5*\cellsize) --++(\cellsize,0);}\missingcell}
\newcommand\vcell{\cellextra{\draw (-0.5*\cellsize,-0.5*\cellsize) --++(0,\cellsize);}\missingcell}
\newcommand\vhcell{\cellextra{\draw (-0.5*\cellsize,-0.5*\cellsize) --++(0,\cellsize) --++(\cellsize,0);}\missingcell}
\def\nwstripecolor{}
\def\nestripecolor{}
\def\stripewidth{0.4pt}
\def\stripesep{0.2}
\newcommand\nwstripedcell{\cellextra{
\path;
\begin{scope}[shift={(-0.5*\cellsize,-0.5*\cellsize)}]
\clip (0,0) rectangle ++(\cellsize,\cellsize);
\foreach\i in {0,\stripesep,...,2} \draw[\nwstripecolor,line width=\stripewidth] (\i*\cellsize,0) -- (0,\i*\cellsize);
\end{scope}
}}
\newcommand\nestripedcell{\cellextra{
\path;
\begin{scope}[shift={(-0.5*\cellsize,0.5*\cellsize)}]
\clip (0,0) rectangle ++(\cellsize,-\cellsize);
\foreach\i in {0,\stripesep,...,2} \draw[\nestripecolor,line width=\stripewidth] (\i*\cellsize,0) -- (0,-\i*\cellsize);
\end{scope}
}}
\newcommand\doublestripedcell{\cellextra{
\path;
\begin{scope}[shift={(-0.5*\cellsize,-0.5*\cellsize)}]
\clip (0,0) rectangle ++(\cellsize,\cellsize);
\foreach\i in {0,\stripesep,...,2} \draw[\nwstripecolor,line width=\stripewidth] (\i*\cellsize,0) -- (0,\i*\cellsize);
\end{scope}
\begin{scope}[shift={(-0.5*\cellsize,0.5*\cellsize)}]
\clip (0,0) rectangle ++(\cellsize,-\cellsize);
\foreach\i in {0,\stripesep,...,2} \draw[\nestripecolor,line width=\stripewidth] (\i*\cellsize,0) -- (0,-\i*\cellsize);
\end{scope}
}}

\def\activate#1{\begingroup
  \lccode`\~=`#1%
  \lowercase{\endgroup \let~#1}%
  \catcode`#1=13\relax}
\activate &
\newcommand\tableau[1]{\tikz[baseline=0]
\node[tableau]{#1};}

\[
\tikz[every matrix/.style={nodes in empty cells=false}]{
\tableau{4&&&\\3&4\\2&3&4\\1&1&2&2\\}
}
\]
where the $i^{\rm th}$ row of the pattern is the shape of the boxes with labels $\le i$.

In general, applying this mapping to strict Gelfand--Tsetlin patterns results in SSYTs with triangular shape
and one additional inequality due to strictness, namely that antidiagonals are weakly increasing: removing the boxes for
clarity, one gets a new triangular array of integers of the form
\[
\begin{tikzpicture}
\useasboundingbox (-3,.4) rectangle (3,-4.1);
\matrix[matrix of math nodes,row sep={-1cm,between origins},column sep={.8cm,between origins}] {
1 &\le& b_{n-1,1} &\le& \cdots &\le& b_{1,1}
\\
2&\le& b_{n-2,2} &\quad \cdots&\quad b_{2,2}
\\
\vdots&&&\ddots
\\
n-1&\le&\quad b_{n-1,n-1}
\\
n
\\
};
\begin{scope}[yscale=-1,yshift=1.8cm]
\node[rotate=90] at (-2.3,1.5) {$<$};
\node[rotate=90] at (-0.7,1.5) {$<$};
\node[rotate=-45] at (-1.5,1.5) {$\ge$};
\node[rotate=-45] at (1.7,1.5) {$\ge$};
\node[rotate=90] at (-2.3,0.5) {$<$};
\node[rotate=-45] at (.6,0.45) {$\ge$};
\node[rotate=90] at (-2.3,-0.55) {$<$};
\node[rotate=-45] at (-.35,-0.45) {$\ge$};
\node[rotate=90] at (-2.3,-1.5) {$<$};
\node[rotate=-45] at (-1.5,-1.5) {$\ge$};
\end{scope}
\end{tikzpicture}
\]
which obeys the same rules (up to a rotation) that define monotone triangles. In fact, the new monotone triangle obtained this way
is nothing but the list of subsets of right-pointing arrows along each column (numbered from bottom to top).


\subsubsection{Fully Packed Loops}\label{ssec:fpl}
Another bijection of a fairly different nature is, given a DWBC configuration, to replace the two
types of decorations of the edges (arrow pointing one way or the other) with another type, which we shall
depict with two colours, say red or blue.
We require that at each vertex, arrows pointing in the same direction (in or out)
should be mapped to the same coloured state. By looking at the edge between two neighbouring vertices,
we conclude that the correspondence should be different depending on whether a vertex is in the odd or even
sublattice of the square lattice; indicating the sublattices by colouring vertices red or blue alternatingly,
the final mapping is that edges acquire the colour of the vertex that they point to, e.g.,
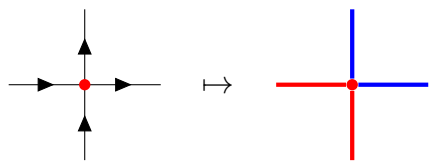
\[
\tikz[baseline=-3pt]{
\draw[sixvin] (-1,0) -- (0,0); \draw[sixvin] (0,-1) -- (0,0); \draw[sixvout] (1,0) -- (0,0); \draw[sixvout] (0,1) -- (0,0);
\node[circle,fill=red,inner sep=1.5pt] at (0,0) {};
}
\quad\mapsto\quad
\tikz[baseline=-3pt]{
\node[circle,fill=red,inner sep=1.5pt] (O) at (0,0) {};
\begin{scope}[ultra thick]
\draw[red] (O) -- (0,-1);
\draw[red] (O) -- (-1,0);
\draw[blue] (O) -- (0,1);
\draw[blue] (O) -- (1,0);
\end{scope}
}
\]

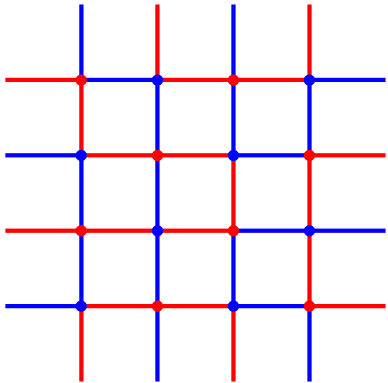
\begin{figure}
\begin{tikzpicture}
\begin{scope}[ultra thick]
\draw[blue] (0,1) -- (1,1);
\draw[red] (1,1) -- (2,1);
\draw[red] (5,1) -- (4,1);
\draw[blue] (4,1) -- (3,1);
\draw[red] (3,1) -- (2,1);
\draw[red] (0,2) -- (1,2);
\draw[red] (2,2) -- (1,2);
\draw[red] (2,2) -- (3,2);
\draw[blue] (3,2) -- (4,2);
\draw[blue] (5,2) -- (4,2);
\draw[blue] (0,3) -- (1,3);
\draw[red] (1,3) -- (2,3);
\draw[red] (3,3) -- (2,3);
\draw[blue] (4,3) -- (3,3);
\draw[red] (5,3) -- (4,3);
\draw[red] (0,4) -- (1,4);
\draw[blue] (1,4) -- (2,4);
\draw[red] (2,4) -- (3,4);
\draw[red] (4,4) -- (3,4);
\draw[blue] (5,4) -- (4,4);
\draw[red] (1,1) -- (1,0);
\draw[blue] (1,2) -- (1,1);
\draw[blue] (1,2) -- (1,3);
\draw[red] (1,3) -- (1,4);
\draw[blue] (1,4) -- (1,5);
\draw[blue] (2,1) -- (2,0);
\draw[blue] (2,1) -- (2,2);
\draw[blue] (2,3) -- (2,2);
\draw[blue] (2,3) -- (2,4);
\draw[red] (2,4) -- (2,5);
\draw[red] (3,1) -- (3,0);
\draw[blue] (3,2) -- (3,1);
\draw[red] (3,3) -- (3,2);
\draw[blue] (3,4) -- (3,3);
\draw[blue] (3,4) -- (3,5);
\draw[blue] (4,1) -- (4,0);
\draw[red] (4,2) -- (4,1);
\draw[red] (4,2) -- (4,3);
\draw[blue] (4,3) -- (4,4);
\draw[red] (4,4) -- (4,5);
\end{scope}
\matrix [matrix of nodes,nodes in empty cells,row sep={1cm,between origins},column sep={1cm,between origins},
every node/.style={circle,fill,inner sep=1.5pt,
\ifodd\pgfmatrixcurrentcolumn 
\ifodd\pgfmatrixcurrentrow red\else blue\fi
\else
\ifodd\pgfmatrixcurrentrow blue\else red\fi
\fi
}
] at (2.5,2.5){
&&&
\\
&&&
\\
&&&
\\
&&&
\\
};
\end{tikzpicture}
\caption{An example of Fully Packed Loop configuration.}
\label{fig:fpl}
\end{figure}

Such configurations are called {\em Fully Packed loops}\/ (FPLs).
The result on our running example is shown on Fig.~\ref{fig:fpl}.
The arrow conservation rule
becomes the rule that every vertex must have two adjacent edges of each colour.
The boundary conditions are alternatingly blue and red.

Note that red (resp.\ blue) lines form uninterrupted paths going from one boundary edge to another
(plus possibly closed loops); this leads to the possibility of refined counting beyond what we have
considered so far. We shall discuss this further in \S\ref{sec:RS}.

\subsubsection{Bumpless pipe dreams}\label{ssec:bpd}
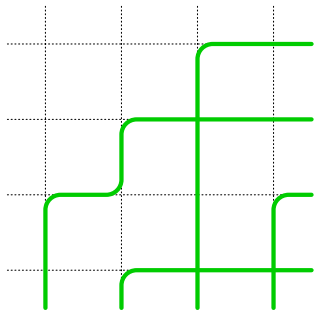
\begin{figure}
\begin{tikzpicture}[line cap=round]
\draw[empty] (0.5,0.5) grid (4.5,4.5);
\begin{scope}[ultra thick,mygreen,rounded corners=2mm]
\draw (1,0.5) -- (1,2) -- (2,2) -- (2,3) -- (4.5,3) ;
\draw (2,0.5) -- (2,1) -- (4.5,1) ;
\draw (3,0.5) -- (3,4) -- (4.5,4) ;
\draw (4,0.5) -- (4,2) -- (4.5,2) ;
\end{scope}
\end{tikzpicture}
\caption{An example of a Bumpless Pipe Dream.}
\label{fig:bpd}
\end{figure}

In a similar vein, let us revisit the lattice path representation of \S\ref{ssec:lp}:
we now declare that paths {\em cross}\/ at vertices rather than osculate. To conform with the
conventions of the literature, we also switch occupied and empty edges. The result is depicted
on Figure~\ref{fig:bpd} for our running example; in general, one obtains a bijection between
$\text{DWBC}_n$ and what is known as Bumpless Pipe Dreams \cite{LLS-BPD,Wei-BPD}.

Green paths now go from the East boundary to the South boundary, giving rise to a permutation
(numbering rows from top to bottom and columns from left to right); in the example, it is $3142$.
(If the configuration is a ``rook placement'', this is the
natural permutation that we have already associated to it in \S\ref{ssec:rook},
so every permutation appears at least once). 
This means that we again have a refined counting by permutation;
we shall not develop this further and refer to the literature \cite{Las-CY,LLS-BPD,Wei-BPD}.

\subsubsection{Domino Tilings}\label{ssec:dt}
\begin{figure}
\begin{tikzpicture}
\begin{scope}[thick]
\draw (1,0) -- (1.5,.5) -- (2,0) -- (2.5,.5) -- (3,0) -- (3.5,.5) -- (4,0) -- (5,1) -- (4.5,1.5) -- (5,2) -- (4.5,2.5) -- (5,3) -- (4.5,3.5) -- (5,4) -- (4,5) -- (3.5,4.5) -- (3,5) -- (2.5,4.5) -- (2,5) -- (1.5,4.5) -- (1,5) -- (0,4) -- (.5,3.5) -- (0,3) -- (.5,2.5) -- (0,2) -- (.5,1.5) -- (0,1) -- cycle;
\draw (1.5,.5) -- (.5,1.5) (3.5,.5) -- (4.5,1.5) (1.5,4.5) -- (.5,3.5) (3.5,4.5) -- (4.5,3.5) (.5,1.5) -- (1.5,2.5) (1.5,1.5) -- (2.5,.5) (1,1) -- (4,4) (1,3) -- (3,1) (.5,2.5) -- (2.5,4.5) (1,4) -- (4,1) (2.5,.5) -- (3.5,1.5) (4.5,1.5) -- (3.5,2.5) (3,2) -- (4,3) (4.5,2.5) -- (3.5,3.5) (2,4) -- (3,3) (2.5,3.5) -- (3.5,4.5);
\end{scope}
\draw[dotted,xshift=.5cm,yshift=.5cm] (0,0) grid (4,4);
\end{tikzpicture}
\caption{An example of Domino Tiling of the Aztec Diamond.}\label{fig:dt}
\end{figure}
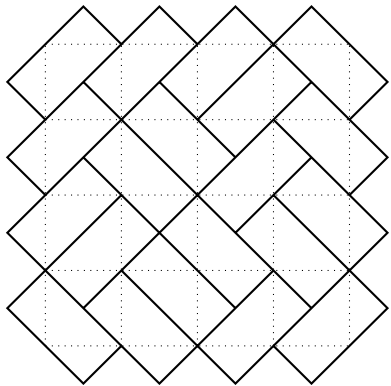

Finally, a mapping that is {\em not}\/ bijective is the following. A {\em Domino Tiling of the Aztec Diamond} of size $n$
is a filling of a staircase domain as on Fig.~\ref{fig:dt} with dominos, i.e., $1\times 2$ rectangles; denote their set $\text{DT}_n$. This terminology is due to Propp, who studied extensively domino tilings of the Aztec diamond \cite{EKLP1,EKLP2,JPS}, in particular because they were the first known model to exhibit the ``limiting shape'' phenomenon as $n\to\infty$. Note that such domino tilings had already appeared in the physics literature under the guise of {\em dimer model},
cf \cite{GCZ} which conjectured
$\#\text{DT}_n$ and pointed out the sensitivity of the model to varying boundary conditions.


The dotted grid is a helper for the mapping which we define now.
Given a tiling of $\text{DT}_n$, apply the following local substitutions:
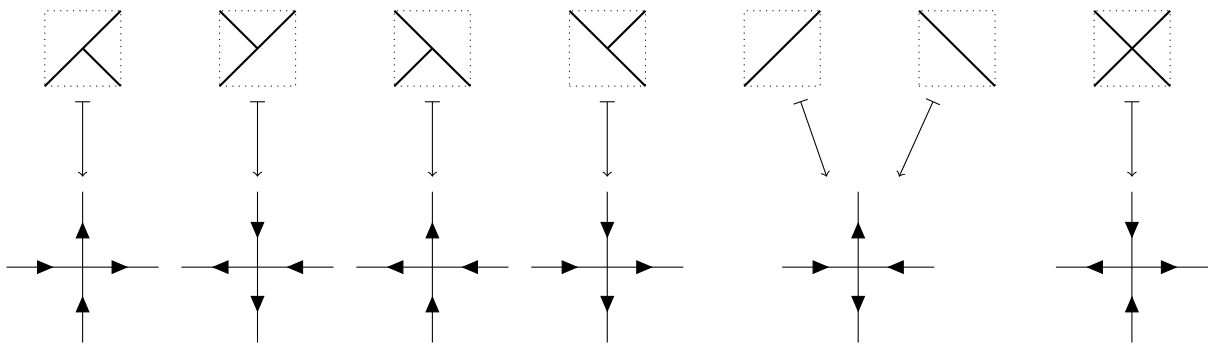
\[
\begin{tikzpicture}[every node/.style={outer sep=2mm}]
\matrix[column sep={2.3cm,between origins}] {
\node[rectangle,dotted,draw,minimum size=1cm] (a1) {};
\draw[thick] (-.5,-.5) -- (.5,.5) (0,0) -- (.5,-.5);
&
\node[rectangle,dotted,draw,minimum size=1cm] (a2) {};
\draw[thick] (-.5,-.5) -- (.5,.5) (0,0) -- (-.5,.5);
&
\node[rectangle,dotted,draw,minimum size=1cm] (b1) {};
\draw[thick] (-.5,.5) -- (.5,-.5) (0,0) -- (-.5,-.5);
&
\node[rectangle,dotted,draw,minimum size=1cm] (b2) {};
\draw[thick] (-.5,.5) -- (.5,-.5) (0,0) -- (.5,.5);
&
\node[rectangle,dotted,draw,minimum size=1cm] (c1a) {};
\draw[thick] (-.5,-.5) -- (.5,.5);
&
\node[rectangle,dotted,draw,minimum size=1cm] (c1b) {};
\draw[thick] (-.5,.5) -- (.5,-.5);
&
\node[rectangle,dotted,draw,minimum size=1cm] (c2) {};
\draw[thick] (-.5,-.5) -- (.5,.5);
\draw[thick] (-.5,.5) -- (.5,-.5);
\\
};
\node[below=of a1,inner sep=1cm] (a1') {};
\draw[sixvout] (a1'.center) -- ++(0,-1);
\draw[sixvout] (a1'.center) -- ++(-1,0);
\draw[sixvin] (a1'.center) -- ++(0,1);
\draw[sixvin] (a1'.center) -- ++(1,0);
\draw[|->] (a1) -- (a1');
\node[below=of a2,inner sep=1cm] (a2') {};
\draw[sixvin] (a2'.center) -- ++(0,-1);
\draw[sixvin] (a2'.center) -- ++(-1,0);
\draw[sixvout] (a2'.center) -- ++(0,1);
\draw[sixvout] (a2'.center) -- ++(1,0);
\draw[|->] (a2) -- (a2');
\node[below=of b1,inner sep=1cm] (b1') {};
\draw[sixvout] (b1'.center) -- ++(0,-1);
\draw[sixvin] (b1'.center) -- ++(-1,0);
\draw[sixvin] (b1'.center) -- ++(0,1);
\draw[sixvout] (b1'.center) -- ++(1,0);
\draw[|->] (b1) -- (b1');
\node[below=of b2,inner sep=1cm] (b2') {};
\draw[sixvin] (b2'.center) -- ++(0,-1);
\draw[sixvout] (b2'.center) -- ++(-1,0);
\draw[sixvout] (b2'.center) -- ++(0,1);
\draw[sixvin] (b2'.center) -- ++(1,0);
\draw[|->] (b2) -- (b2');
\node[below=of c1a,inner sep=1cm,xshift=1cm] (c1') {};
\draw[sixvin] (c1'.center) -- ++(0,-1);
\draw[sixvout] (c1'.center) -- ++(-1,0);
\draw[sixvin] (c1'.center) -- ++(0,1);
\draw[sixvout] (c1'.center) -- ++(1,0);
\draw[|->] (c1a) -- (c1');
\draw[|->] (c1b) -- (c1');
\node[below=of c2,inner sep=1cm] (c2') {};
\draw[sixvout] (c2'.center) -- ++(0,-1);
\draw[sixvin] (c2'.center) -- ++(-1,0);
\draw[sixvout] (c2'.center) -- ++(0,1);
\draw[sixvin] (c2'.center) -- ++(1,0);
\draw[|->] (c2) -- (c2');
\end{tikzpicture}
\]
The result is a configuration in $\text{DWBC}_n$.
Applying this mapping to Fig.~\ref{fig:dt} takes us back to our running example of Fig.~\ref{fig:6vdwbc}.

This mapping is clearly not bijective since two local configurations are sent to one. From the point of view
of partition functions, this can be absorbed in a doubling of the appropriate Boltzmann weight, namely
$c_1$. In particular the number of domino tilings of the Aztec Diamond is equal to the DWBC partition function
with Boltzmann weights $a=b=c_2=1$, $c_1=2$.

\subsection{The Izergin determinant formula}
We now proceed to compute the partition function of the six-vertex model with DWBC.

We choose the following convenient parameterisation of the local Boltzmann weights:
\begin{alignat}{2}
\notag
a(x,y)&=q x - q^{-1}y
\\\label{eq:6vweights}
b(x,y)&=x-y
\\\notag
c_1(x,y)&=(q-q^{-1})y&\qquad c_2(x,y)&=(q-q^{-1})x
\end{alignat}
(i.e., $c(x,y)=(q-q^{-1})(xy)^{1/2}$, but we have avoided the use of square roots by introducing $c_1$, $c_2$).
Here $q$, $x$ and $y$ are formal parameters, which can be chosen to be nonzero complex numbers.
If we want the Boltzmann weights to be real (up to overall normalisation), then one should choose $q$ and $z$
to be either real or of modulus one. This traditionally leads to a division into 3 regimes ($q$
real positive, real negative,
of modulus one) which will not be discussed here. \rem{ref to other chapter?}

Though all three weights depend on all parameters $q,x,y$, we emphasise in the notation the dependence
on $x,y$; the reason is that
it is useful to make the model {\em inhomogeneous}\/ by varying $x$ and $y$ (but not $q$) depending on the
row/column of the vertex. Namely, at the vertex at row $i$ and column $j$ of the square $n\times n$ domain,
we use the Boltzmann weights \eqref{eq:6vweights} with
the substitution $x=x_i$, $y=y_j$, where $x_1,\ldots,x_n,y_1,\ldots,y_n$ are fixed parameters.

With this choice of Boltzmann weights, denoting $Z_n(x_1,\ldots,x_n;y_1,\ldots,y_n)$ the partition function
of the six-vertex model with DWBC:
\[
Z_n(x_1,\ldots,x_n;y_1,\ldots,y_n)
=
\sum_{\mathcal C \in \text{DWBC}_n} \prod_{i,j=1}^n
\left\{\begin{smallmatrix}a\\b\\c_1\\c_2\\
  \end{smallmatrix}
  \right\}(x_i,y_j)
\]
 one has the following beautiful formula:
\begin{theorem}[Izergin \cite{Iz-6V}]\label{thm:DWBC}
\[
Z_n(x_1,\ldots,x_n;y_1,\ldots,y_n)
=
\prod_{j=1}^n y_j\,
\frac{\prod_{i,j=1}^n (x_i-y_j)(qx_i-q^{-1}y_j)}{\prod_{1\le i<j\le n} (x_i-x_j)(y_j-y_i)}
\det_{i,j=1,\ldots,n} \left(\frac{q-q^{-1}}{(x_i-y_j)(qx_i-q^{-1}y_j)}\right)
\]
\end{theorem}

The proof is by now standard and can be found in e.g.~\cite{ICK,artic12}, and we only give a sketch of it.
The main ingredient is that the local Boltzmann weights satisfy the {\em Yang--Baxter equation} (YBE),
which using our parameterisation takes the simple form:
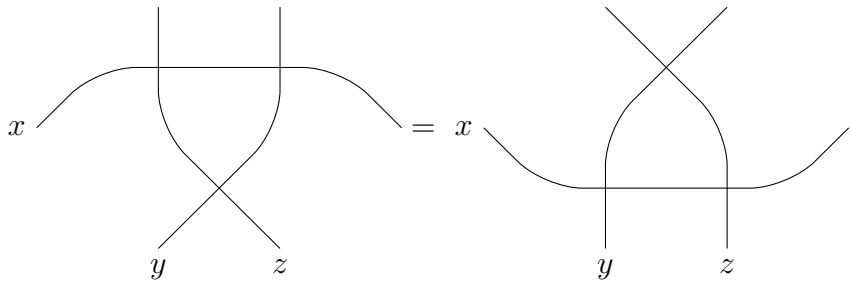
\[
\begin{tikzpicture}[baseline=-3pt,rotate=90,rounded corners=5mm,scale=.8]
\draw (-2,-1) node[below] {$z$} -- (0,1) -- (2,1);
\draw (-2,1) node[below] {$y$} -- (0,-1) -- (2,-1);
\draw (0,-3) -- (1,-2) -- (1,2) -- (0,3) node[left] {$x$};
\end{tikzpicture}
=
\begin{tikzpicture}[baseline=-3pt,rotate=90,rounded corners=5mm,scale=.8]
\draw (-2,-1) node[below] {$z$} -- (0,-1) -- (2,1);
\draw (-2,1) node[below] {$y$} -- (0,1) -- (2,-1);
\draw (0,-3) -- (-1,-2) -- (-1,2) -- (0,3) node[left] {$x$};
\end{tikzpicture}
\]
Each picture stands for the corresponding partition function; more precisely,
the convention is that the labelling of external edges is arbitrary but fixed (so the YBE is really
$2^6=64$ equations, though many are trivial or redundant), whereas the arrows
of the internal edges are summed over. The Boltzmann weights used at each vertex involve
the parameters attached to each line, namely, $(x,y)$, $(x,z)$, $(y,z)$. 

The YBE can be used repeatedly, showing symmetry of $Z_n$ under interchange of the $y_j$s (and similarly
for the $x_i$s):
\begin{multline*}
a(y_{j+1},y_j)Z_n(\ldots,y_j,y_{j+1},\ldots)
=
\\
\begin{tikzpicture}[baseline=(current  bounding  box.center),rounded corners]
\foreach\x in {1,4} \draw[sixvout=.12,sixvin=.92] (\x,0) -- (\x,5);
\foreach\y in {1,2,3,4} \draw[sixvin=.12,sixvout=.92] (0,\y) -- (5,\y);
\draw[sixvout=.06,sixvin=.94] (2,-1) node[below] {$y_{j+1}$} -- (3,0) -- (3,5);
\draw[sixvout=.06,sixvin=.94] (3,-1) node[below] {$y_j$} -- (2,0) -- (2,5);
\end{tikzpicture}
= 
\begin{tikzpicture}[baseline=(current  bounding  box.center),rounded corners]
\foreach\x in {1,4} \draw[sixvout=.12,sixvin=.92] (\x,0) -- (\x,5);
\foreach\y in {1,2,3,4} \draw[sixvin=.12,sixvout=.92] (0,\y) -- (5,\y);
\draw[sixvout=.09,sixvin=.96] (2,0) node[below] {$y_{j+1}$} -- (2,5) -- (3,6);
\draw[sixvout=.09,sixvin=.96] (3,0) node[below] {$y_j$} -- (3,5) -- (2,6);
\end{tikzpicture}
\\
= a(y_{j+1},y_j)Z_n(\ldots,y_{j+1},y_{i},\ldots)
\end{multline*}
Finally, one finds by inspection that setting $x_n=y_n$ leads to a recurrence relation for $Z_n$, namely
\[
Z_n(x_1,\ldots,x_n;y_1,\ldots,y_n)|_{x_n=y_n} = c_1(y_n,y_n)\prod_{i=1}^{n-1} a(x_i,y_n) \prod_{j=1}^{n-1} a(y_n,y_j)\ 
Z_{n-1}(x_1,\ldots,x_{n-1};y_1,\ldots,y_{n-1})
\]
By symmetry, this provides the value of $Z_n$ at $n$ distinct specialisations of $x_n$, and since it is not
hard to show that $Z_n$ is a polynomial of degree at most $n-1$ in $x_n$ (because each row has
at least one vertex of type $c_1$), this specifies it uniquely. Only remains to check
that the provided expression satisfies this recurrence relation.

\subsection{Homogeneous limit}
Note that the expression in Theorem~\ref{thm:DWBC} is indeterminate when all $x_i$s (or $y_i$s) are equal.
By applying L'H\^opital's rule, one can find an expression for $Z_n$ when $x_i\to x$ and $y_i\to y$,
as an $n\times n$ determinant.

This determinant can be re-expressed in various ways, see e.g. \cite{artic13,CP-6v,artic55}.
It can then be evaluated at special values of $a,b,c$, cf \cite{Kup-ASM,CP-orthog1};
two such evaluations will be discussed below.

Finally, one can compute asymptotics of $Z_n$ in the ``thermodynamic limit'' $n\to\infty$; we shall not discuss
this here and refer to the extensive literature \cite{artic12,artic13,CP2,Bleher-review}.

\subsubsection{The ice point and exact enumeration}
By setting $a=b=c=1$, the partition function is simply the cardinality of $\text{DWBC}_n$, or of any of the sets
in bijection such as $\text{ASM}_n$. This can be achieved by choosing parameters $q=e^{i\pi/3}$,
$x=e^{-i\pi/6}/\sqrt3$, $y=-qx$ (or their complex conjugate).

Following \cite{Strog-IK,Oka}, we may compute
$Z_n$ as follows: we first set $q=e^{i\pi/3}$, leaving the parameters $x_i$ and $y_j$ free.
Then it turns out that the recurrence relation becomes identical to one satisfied by the {\em Schur polynomial}
$s_{\lambda^{(n)}}$
associated to the partition $\lambda^{(n)}=(n-1,n-1,n-2,n-2,\ldots,1,1)$. We refer to \S\ref{ssec:schur} below for
an introduction to Schur polynomials. Taking care of various prefactors in the recurrence relations, we
have the identification
\[
Z_n(x_1,\ldots,x_n;y_1,\ldots,y_n)|_{q=e^{i\pi/3}}=q^{-n(n-1)/2}(q-q^{-1})^{n} s_{\lambda^{(n)}}(-qx_1,\ldots,-qx_n,y_1,\ldots,y_n)
\]
Now we can specialise $x_i=x=e^{-\pi i/6}/\sqrt3$, $y_j=-qx$ and find, using the homogeneity property of Schur polynomials,
$Z_n=(-q)^{-n(n-1)/2}(q-q^{-1})^{n} (-e^{i\pi/6}/\sqrt{3})^{n(n-1)}s_{\lambda^{(n)}}(1,\ldots,1)$.
Simplifying the prefactor and taking into account the extra phase arising from $(c_1/c_2)^{n/2}$, we obtain
the formula
\[
\# \text{DWBC}_n = \#\text{ASM}_n = 3^{-n(n-1)/2} s_{\lambda^{(n)}}(\underbrace{1,\ldots,1}_{2n})
\]
The evaluation of a Schur polynomial at $1,\ldots,1$ is known (see \eqref{eq:dimschur});
simplifying the resulting product, we obtain our final identity
\[
\#\text{ASM}_n
=\prod_{i=0}^{n-1} \frac{(3i+1)!}{(n+i)!} = 1,2,7,42,429\ldots
\]
which is the famous formula for the number of ASMs (OEIS sequence \href{https://oeis.org/A005130}{A005130}).

\subsubsection{The free fermion point}\label{ssec:ff}
Consider the $\lambda$-determinant of \eqref{eq:ldet2}. In order to achieve the Boltzmann weights
there, one can pick $x=(1-i\sqrt{\lambda})/2$, $y=-(1+i\sqrt{\lambda})/2$, $q=i$, resulting
in $a=\sqrt{\lambda}$, $b=1$, $c=\sqrt{c_1c_2}=\sqrt{1+\lambda}$. 

As a special case, if we set $\lambda=1$, we recognize
the choice of weights $a=b=1$, $c=\sqrt 2$
that is relevant to the counting of domino tilings, cf \S\ref{ssec:dt}.


As in the previous section, 
let us first set $q$ to its value, namely $i$ (this is known as the {\em free fermion point}\/ of the six-vertex model;
some justification for this terminology will be given in \S\ref{ssec:pDWBC}).
The determinant in Theorem~\ref{thm:DWBC}
then turns into a Cauchy determinant, and the partition function factorises as
\begin{equation}\label{eq:ff}
Z_n(x_1,\ldots,x_n;y_1,\ldots,y_n)|_{q=i}
=
(2i)^n (-1)^{n(n-1)/2} \prod_{j=1}^n y_j \prod_{1\le i<j\le n} (x_i+x_j)(y_i+y_j)
\end{equation}
We can now take the homogeneous limit as above; reintroducing
the extra factor $(-i+\sqrt\lambda)^{-n}$ needed for the separate weights $c_1$ and $c_2$ as in \eqref{eq:ldet2}, we obtain
\[
\ldet I_n = (1+\lambda)^{n(n-1)/2}
\]
an equality which can be easily derived from the definition of the $\lambda$-determinant.

To obtain the weights of \S\ref{ssec:dt}, we need to set $\lambda=1$ and to multiply by
the factor $(1+i)^n$, resulting in
\[
\# \text{DT}_n = 2^{n(n+1)/2}
\]
which is the well-known formula for the number of Domino Tilings of the Aztec Diamond \cite{EKLP1}.

\subsection{The Razumov--Stroganov correspondence}\label{sec:RS}
Let us now redraw the $7$ configurations in $\text{DWBC}_3$ as Fully Packed Loop configurations, cf \S\ref{ssec:fpl}:
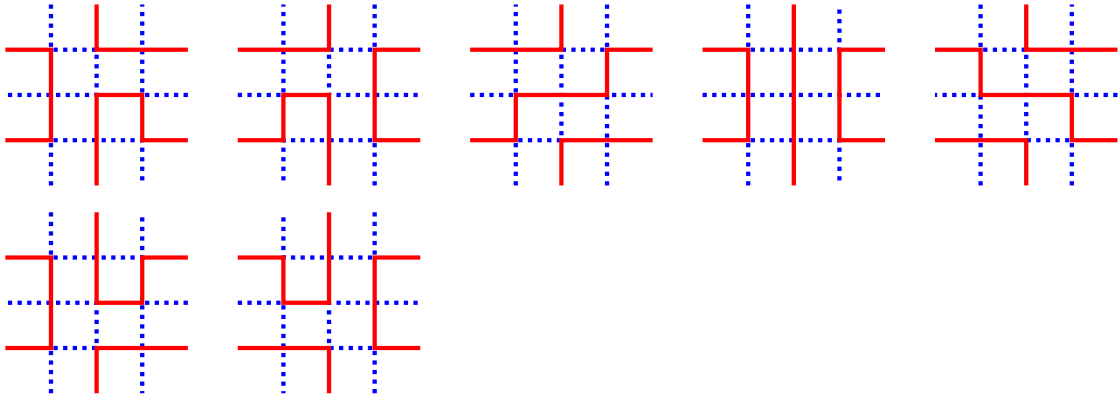
\begin{center}
\begin{tikzpicture}[ultra thick]
\matrix[every cell/.style={scale=.6},column sep=.6cm,row sep=.3cm]{
\draw[dotted,blue] (1,0) -- (1,1) -- (3,1) -- (3,0) (1,4) -- (1,3) -- (2,3) -- (2,2) -- (0,2) (4,2) -- (3,2) -- (3,4);
\draw[red] (0,1) -- (1,1) -- (1,3) -- (0,3) (4,1) -- (3,1) -- (3,2) -- (2,2) -- (2,0) (2,4) -- (2,3) -- (4,3);
&
\begin{scope}[xscale=-1]
\draw[dotted,blue] (1,0) -- (1,1) -- (3,1) -- (3,0) (1,4) -- (1,3) -- (2,3) -- (2,2) -- (0,2) (4,2) -- (3,2) -- (3,4);
\draw[red] (0,1) -- (1,1) -- (1,3) -- (0,3) (4,1) -- (3,1) -- (3,2) -- (2,2) -- (2,0) (2,4) -- (2,3) -- (4,3);
\end{scope}
&
\draw[dotted,blue] (1,0) -- (1,1) -- (2,1) -- (2,3) -- (3,3) -- (3,4) (3,0) -- (3,2) -- (4,2) (0,2) -- (1,2) -- (1,4);
\draw[red] (0,1) -- (1,1) -- (1,2) -- (3,2) -- (3,3) -- (4,3) (0,3) -- (2,3) -- (2,4) (2,0) -- (2,1) -- (4,1);
&
\draw[dotted,blue] (1,0) -- (1,1) -- (3,1) -- (3,0) (1,4) -- (1,3) -- (3,3) -- (3,4) (0,2) -- (4,2);
\draw[red] (0,1) -- (1,1) -- (1,3) -- (0,3) (4,1) -- (3,1) -- (3,3) -- (4,3) (2,0) -- (2,4);
&
\begin{scope}[xscale=-1]
\draw[dotted,blue] (1,0) -- (1,1) -- (2,1) -- (2,3) -- (3,3) -- (3,4) (3,0) -- (3,2) -- (4,2) (0,2) -- (1,2) -- (1,4);
\draw[red] (0,1) -- (1,1) -- (1,2) -- (3,2) -- (3,3) -- (4,3) (0,3) -- (2,3) -- (2,4) (2,0) -- (2,1) -- (4,1);
\end{scope}
\\
\begin{scope}[yscale=-1]
\draw[dotted,blue] (1,0) -- (1,1) -- (3,1) -- (3,0) (1,4) -- (1,3) -- (2,3) -- (2,2) -- (0,2) (4,2) -- (3,2) -- (3,4);
\draw[red] (0,1) -- (1,1) -- (1,3) -- (0,3) (4,1) -- (3,1) -- (3,2) -- (2,2) -- (2,0) (2,4) -- (2,3) -- (4,3);
\end{scope}
&
\begin{scope}[scale=-1]
\draw[dotted,blue] (1,0) -- (1,1) -- (3,1) -- (3,0) (1,4) -- (1,3) -- (2,3) -- (2,2) -- (0,2) (4,2) -- (3,2) -- (3,4);
\draw[red] (0,1) -- (1,1) -- (1,3) -- (0,3) (4,1) -- (3,1) -- (3,2) -- (2,2) -- (2,0) (2,4) -- (2,3) -- (4,3);
\end{scope}
\\
};
\end{tikzpicture}
\end{center}
Though blue and red paths play symmetric roles, we focus on the red paths, hence the dotting of the blue paths.

We remark that we can group together configurations according to the connectivity of the external (red) edges, i.e., the pairing between them
induced by the paths. In the list above, each column corresponds to a different connectivity. We can then do a refined counting according to
connectivity, leading for $n=3$ to $(2,2,1,1,1)$.

A first observation is that if we rotate the connectivity (as unnatural as it may seem considering we are on the square lattice), then the counting
remains the same: here the $2$s and the $1$s form an orbit each under connectivity rotation.
Such a statement can be proven (for all $n$) \cite{Wieland}, by defining a ``gyration'' operation on FPLs.

Furthermore, Razumov and Stroganov conjectured in \cite{RS-conj} that the vector of FPL numbers, e.g. $(2,2,1,1,1)$ at $n=3$, is the eigenvector of a matrix
that is easy to write explicitly; in fact, this matrix is nothing but the Hamiltonian of a quantum integrable system! See the review \cite{dG-review3}
for a full statement of this correspondence, as well as generalisations.
This conjecture was proven in \cite{CS-RS} by a nontrivial use of the gyration of \cite{Wieland}, though various generalisations are still open.

\section{Symmetric polynomials and quantum integrability}\label{sec:sym}
Let $\Lambda_n=\ZZ[z_1,\ldots,z_n]^{\mathcal S_n}$ be the ring
of {\em symmetric polynomials}\/ in $n$ variables.
The study of $\Lambda_n$ is a classical subject of algebraic combinatorics \cite{Macdonald}.\footnote{A closely related
concept is that of symmetric functions, that is loosely, of symmetric polynomials in an infinite number of variables.
Although such an $n\to\infty$ limit can be taken in our lattice models, we shall not discuss it here.}
The most famous basis of $\Lambda_n$ (as a free $\ZZ$-module) consists of {\em Schur polynomials},
whose definition will be reviewed below.

The first connection to integrability can be traced back to work of the Kyoto school, see
\cite{JM-ff} and references therein,
where it is observed that Schur polynomials are particular solutions of the KP hierarchy. Classical integrability
is not directly related to the present lattice models; however, the key ingredient of this construction is
(quantum mechanical) {\em free fermions}, suggesting that perhaps our free-fermionic six-vertex model
(cf \S\ref{ssec:ff}) may be relevant.

Such a connection between Schur polynomials and a particular case of the free-fermionic six-vertex model,
namely the free fermionic {\em five-vertex} model, was indeed found in \cite{artic46}. The connection to
the work \cite{JM-ff} can be found in \cite{hdr}.

\rem{references will be messy. things seem to have been rediscovered so many times, cf this silly paper
\url{https://arxiv.org/pdf/0910.5288.pdf} and its refs}

The most general free-fermionic six-vertex, as well as a comprehensive overview of the connection to Schur polynomials, can be found in \cite{Nap-Schur}.

It is important to stress that the connection between solvable lattice models and symmetric polynomials
extends beyond free-fermionic models, and there is by now an extensive literature on the subject;
see for example the lattice models of \cite{Motegi-Groth,artic68} related
to symmetric Grothendieck polynomials or of \cite{tsi,artic65} for Hall--Littlewood polynomials; and the two applications
that we mention in what follows, namely the Cauchy identity (\S\ref{ssec:cauchy}) and product rule (\S\ref{ssec:puzzle})
can be extended to these other families of polynomials.

There is also a deep connection to geometry which we cannot develop here; see e.g.~\cite{GK-K,Korff-qC,Korff-GW,artic71}
for a connection between various cohomology rings, their Schubert classes and their polynomial representatives,
to quantum integrable systems. In this language, the example of Schur polynomials
which we develop now is related to the (ordinary) cohomology of Grassmannians.

\subsection{Schur polynomials}\label{ssec:schur}
Define Schur polynomials as follows. A {\em partition}\/ with $n$ parts is a weakly decreasing
sequence of nonnegative integers $\lambda=(\lambda_1\ge\cdots\ge\lambda_n\ge0)$ (zero parts can be omitted, so that partitions with less than
$n$ parts are implicitly padded with zeroes). To such a $\lambda$ we associate the polynomial
\begin{equation}\label{eq:JW}
s_\lambda(z_1,\ldots,z_n)=\frac{\det_{i,j=1,\ldots,n} (z_i^{\lambda_j+n-j})}{\prod_{1\le i<j\le n}(z_i-z_j)}
\end{equation}
(Jacobi's bilalternant formula, which is a special case of Weyl's character formula).

It is easy to see that $s_\lambda$ is a homogeneous symmetric polynomial in the $z_i$, of degree $|\lambda|=\sum_{i=1}^n \lambda_i$.
Also note the specialisation
\begin{equation}\label{eq:dimschur}
s_{\lambda}(1,\ldots,1) = \prod_{1\le i<j\le n} \frac{\lambda_i-i-\lambda_j+j}{j-i}
\end{equation}
(hook-length formula).

We shall take \eqref{eq:JW} as a definition of $s_\lambda$ and obtain as a byproduct of our lattice model formulation of $s_\lambda$ two equivalent expressions.

The first one is the Jacobi--Trudi formula. Introduce the generating series
\begin{equation}\label{eq:h0}
h(u)=\prod_{i=1}^n \frac{1}{1-uz_i}\qquad h(u)=\sum_{k=-\infty}^\infty h_k u^k
\end{equation}
Then one has
\begin{equation}\label{eq:JT}
s_\lambda(z_1,\ldots,z_n)=\det_{i,j=1,\ldots,n} h_{\lambda_j-j+i}
\end{equation}

Secondly, recall that a partition can be viewed as a Ferrers diagram, e.g., 
$(4,1)=
\tikz[every matrix/.style={nodes in empty cells=false}]{
\tableau{\ &&&&\\\ &\ &\ &\ &\\};
}$.

A {\em Semi-Standard Young tableau}\/ with shape $\lambda$ and alphabet $\{1,\ldots,n\}$ is a filling of the diagram of $\lambda$ with
those integers which is weakly increasing along rows and strictly decreasing along columns. Denote their set
$\text{SSYT}_n(\lambda)$. Then one has
\begin{equation}\label{eq:SSYT}
s_\lambda(z_1,\ldots,z_n)=\sum_{T\in\text{SSYT}_n(\lambda)} \prod_{i=1}^n z_i^{\#\{\text{occurrences of }i\text{ in }T\}}
\end{equation}

\subsection{Partial DWBC}\label{ssec:pDWBC}
We consider the six-vertex model with boundary conditions which generalise the DWBC.
Note that all the alternate representations of \S\ref{ssec:alt} can be adapted to the present setting, and in what follows we shall use
the path representation of \S\ref{ssec:lp} to describe our configurations.

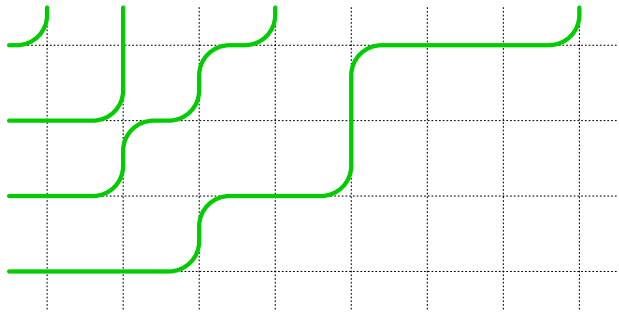
\begin{figure}
\begin{tikzpicture}[line cap=round]
\draw[empty] (0.5,0.5) grid (8.5,4.5);
\begin{scope}[every path/.style={draw=mygreen,ultra thick,rounded corners=4mm}]
\draw (0.5,1) -- (3,1) -- (3,2) -- (5,2) -- (5,4) -- (8,4) -- (8,4.5) ;
\draw (0.5,2) -- (2,2) -- (2,3) -- (3,3) -- (3,4) -- (4,4) -- (4,4.5) ;
\draw (0.5,3) -- (2,3) -- (2,4.5) ;
\draw (0.5,4) -- (1,4) -- (1,4.5) ;
\end{scope}
\end{tikzpicture}
\caption{An example of a lattice path configuration with top boundary given by $\lambda=(4,1,0,0)$.}
\label{fig:lp2}
\end{figure}

We fix a positive integer $n$ and a partition $\lambda_1\ge\cdots\ge\lambda_n\ge0$, as well as an integer $p\ge n+\lambda_1$ (whose value will turn out
to be irrelevant). We then consider the six-vertex model on a $n\times p$ grid, where
the boundary conditions associated to $\lambda$ are the following: on West, South, and East sides, we have the same as for DWBC, that is in terms of paths,
$n$ paths entering from the West and none entering/exiting from South/East; and on the North side, the paths exit in such a way that
the $i^{\rm th}$ path counted from the right is $\lambda_i$ steps to the right of its leftmost position (i.e., exits at column $n+1-i+\lambda_i$).
See Figure~\ref{fig:lp2} for an example.

The Boltzmann weights are the same as in \ref{ssec:ff}, except we use a slightly different normalisation:
\begin{equation}\label{eq:ffwts}
a=1\qquad b=z\qquad c_1=1+z^2\qquad c_2=1
\end{equation}
In terms of paths, this means that vertices where paths go straight through have a weight of $z$, and bends to the right when there is no path touching 
have a weight of $1+z^2$.

We use a superscript $(1)$ to denote partition functions with such a choice of Boltzmann weights, for reasons to become
clear below. Furthermore we want the weights row- (but not column-) dependent:
we set $z=z_i$ at row $i$ of the lattice. We thus define $Z_\lambda^{(1)}(z_1,\ldots,z_n)$ to be the partition function
with boundary conditions associated to $\lambda$. Note that the same YBE argument that we used for DWBC shows here that 
$Z_\lambda^{(1)}(z_1,\ldots,z_n)$ is a {\em symmetric}\/ polynomial in the $z_i$.

\begin{figure}
\begin{tikzpicture}[line cap=round]
\draw[empty] (-3.5,0.5) grid (8.5,4.5);
\draw[dashed] (0.5,0.5) -- (0.5,4.5);
\begin{scope}[every path/.style={draw=mygreen,ultra thick,rounded corners=4mm}]
\draw (0,0.5) -- (0,1) -- (3,1) -- (3,2) -- (5,2) -- (5,4) -- (8,4) -- (8,4.5) ;
\draw (-1,0.5) -- (-1,1) -- (0,1) -- (0,2) -- (2,2) -- (2,3) -- (3,3) -- (3,4) -- (4,4) -- (4,4.5) ;
\draw (-2,0.5) -- (-2,1) -- (-1,1) -- (-1,3) -- (0,3) -- (2,3) -- (2,4.5) ;
\draw (-3,0.5) -- (-3,3) -- (-1,3) -- (-1,4) -- (1,4) -- (1,4.5) ;
\end{scope}
\end{tikzpicture}
\caption{An example of an extended lattice path configuration with top boundary given by $\lambda=(4,1,0,0)$.}
\label{fig:lp3}
\end{figure}
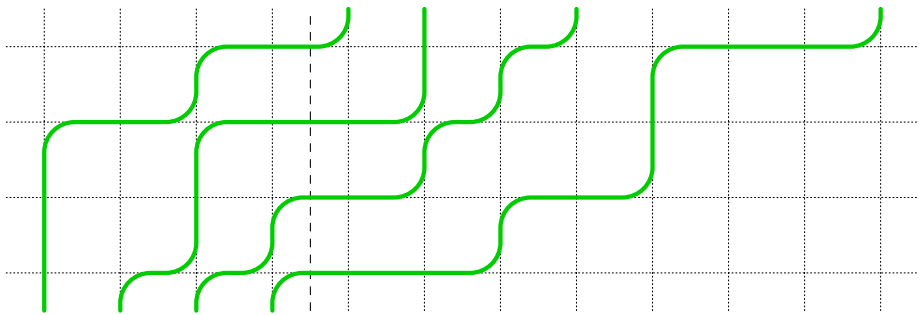

We have the following formula (as a special case of a result first stated in \cite{cslectures}; see also \cite{BBF-Schur} for a similar special case of the latter):
\rem{i'm sure there's a better ref out there}
\begin{theorem}
\begin{equation}\label{eq:pDWBC}
Z^{(1)}_\lambda(z_1,\ldots,z_n)=\prod_{1\le i<j\le n} (1+z_iz_j)\ s_\lambda(z_1,\ldots,z_n)
\end{equation}
\end{theorem}
We give here a sketch of proof. First, it is convenient to extend the configurations: define ``extended'' boundary conditions associated to $\lambda$
by adding $n$ columns to the left, resulting in a $n\times (n+p)$ grid, and having paths start on the South side at the $n$ leftmost locations,
see Figure~\ref{fig:lp3}.

It is easy to see that all paths go through the dashed line separating the new region from the old one, so that
the partition function $\tilde Z^{(1)}_\lambda$ satisfies $\tilde Z^{(1)}_\lambda = Z^{(1)}_n Z^{(1)}_\lambda$,
where $Z^{(1)}_n$ is nothing but the DWBC partition function which we have already computed
in \S\ref{ssec:ff}, cf \eqref{eq:ff}; adapting the result to our new conventions, we find 
$Z^{(1)}_n=\prod_{1\le i\le j\le n} (1+z_iz_j)$, and therefore
\[
\tilde Z^{(1)}_\lambda(z_1,\ldots,z_n) = \prod_{1\le i\le j\le n} (1+z_iz_j) \ Z^{(1)}_\lambda(z_1,\ldots,z_n)
\]
The ``free fermionic'' nature of the Boltzmann weights which we have chosen means that we can apply Wick's theorem,
also known as the LGV formula in the context of lattice paths 
\cite{GV,Lind}: the partition function is the determinant
of partition functions of a single path coming from a given starting point to a given endpoint. Denote by $h^{(1)}_k$
the partition function of a single path whose endpoint is $k+n$ steps to the right of its starting point (because
our Boltzmann weights are independent of the column, only the horizontal distance between the starting and end
points matter; the shift by $n$ is chosen for convenience). Introduce the generating series
$h^{(1)}(u)=\sum_k h^{(1)}_k u^k$. A simple calculation shows that
\[
h^{(1)}(u)=\prod_{i=1}^n \frac{1+u/z_i}{1-uz_i}
\]
Conversely, given $h^{(1)}(u)$, one can recover $h^{(1)}_k$ by writing $h^{(1)}_k = \oint \frac{du}{2\pi i} u^{-k-1} h^{(1)}(u)$ where the contour is around
zero. A standard calculation based on residues then leads to
\begin{align*}
\tilde Z^{(1)}_\lambda(z_1,\ldots,z_n) = \det h^{(1)}_{\lambda_j-j+i}
&=\sum_{\sigma\in\mathcal S_n} (-1)^{\ell(\sigma)}
\oint \prod_{j=1}^n \frac{du_j}{2\pi i} u_j^{-\lambda_j+j-\sigma(j)} \prod_{i,j=1}^n\frac{1+u_j/z_i}{1-u_jz_i} 
\\
&=\oint \prod_{j=1}^n \frac{du_j}{2\pi i} u_j^{-\lambda_j+j} \prod_{1\le j<k\le n}(u_j-u_k) \prod_{i,j=1}^n\frac{1+u_j/z_i}{1-u_jz_i} 
\\
&=\sum_{\sigma\in\mathcal S_n} (-1)^{\ell(\sigma)} z_{\sigma(j)}^{\lambda_j-j} \prod_{1\le j<k\le n}(z_j-z_k) 
\frac{\prod_{i,j=1}^n(1+(z_{\sigma(j)}z_i)^{-1})}{\prod_{\substack{\!\!\!i,j=1\\i\ne\sigma(j)}}^n(1-z_i/z_{\sigma(j)})}
\\
&=\prod_{i,j=1}^n (1+z_iz_j)\, \frac{\det(z_i^{\lambda_j+n-j})}{\prod_{1\le i<j\le n}(z_i-z_j)}
\end{align*}
which is the desired formula by comparing with \eqref{eq:JW}.

\subsection{The five-vertex model limit}
We shall now consider the limit where the $z_i$ are small, so that $Z^{(1)}_\lambda \sim s_\lambda$. In terms of lattice paths,
it means that paths which have as few straight parts as possible are favoured: the leading configurations
are those in which the paths only go East or North-East.

We can formalise this reasoning as follows.
We want to rescale weights \eqref{eq:ffwts} by $z\mapsto \sqrt \alpha z$, but as usual we want to get rid of the square
root by separating $b_1$ and $b_2$. It is easy to see that in a configuration with boundary conditions given by $\lambda$, there are $|\lambda|$ more
vertices of type $b_2$ than of type $b_1$; this means that if we choose the weights
\begin{equation}\label{eq:ffwts2}
a=1\qquad b_1=\alpha z\qquad b_2=z\qquad c_1=1+\alpha z^2\qquad c_2=1
\end{equation}
and denote the corresponding partition functions with the superscript $(\alpha)$, then one has
\[
Z_\lambda^{(\alpha)}(z_1,\ldots,z_n) = \alpha^{-|\lambda|/2} Z_\lambda^{(1)}(\sqrt \alpha z_1,\ldots,\sqrt \alpha z_n)
\]
Plugging this into \eqref{eq:pDWBC} and using the homogeneity property of Schur polynomials,
one finds
\[
  Z^{(\alpha)}_\lambda(z_1,\ldots,z_n)=\prod_{1\le i\le j\le n} (1+\alpha z_iz_j)\ s_\lambda(z_1,\ldots,z_n)
\]
and similarly,
$\tilde Z^{(\alpha)}_\lambda(z_1,\ldots,z_n)=\prod_{i,j=1}^n (1+\alpha z_iz_j)\ s_\lambda(z_1,\ldots,z_n)$.

We can now set $\alpha=0$, which means the vanishing of the vertex $b_1$ (paths cannot go straight North):
the resulting model is called the (free-fermionic)
{\em five-vertex model}.

In this limit, we recover directly the two alternative definitions of Schur polynomials.
Firstly, note that $\tilde Z_\lambda^{(\alpha)}=\det h_{\lambda_j-j+i}^{(\alpha)}$ with the generating series given by
$h^{(\alpha)}(u)=\frac{1+\alpha u/z_i}{1-uz_i}$. By setting $\alpha=0$, we find $s_\lambda = \tilde Z_\lambda^{(0)}
=\det h_{\lambda_j-j+i}^{(0)}$; but $h^{(0)}(u)=h(u)$ (cf \eqref{eq:h0}) and we recover the Jacobi--Trudi formula \eqref{eq:JT}.

\begin{figure}
\begin{tikzpicture}[line cap=round]
\draw[empty] (0.5,0.5) grid (8.5,4.5);
\begin{scope}[every path/.style={draw=mygreen,ultra thick,rounded corners=4mm}]
\draw (0.5,1) -- node[pos=0.2] {1} node[pos=0.6] {1} (3,1) -- (3,2) -- node {2} (5,2) -- (5,3) -- (6,3) -- (6,4) -- node {4} (8,4) -- (8,4.5) ;
\draw (0.5,2) -- node[pos=0.33] {2} (2,2) -- (2,3) -- node {3} (4,3) -- (4,4) -- (5,4) -- (5,4.5) ;
\draw (0.5,3) -- (1,3) -- (1,4) -- node {4} (3,4) -- (3,4.5) ;
\draw (0.5,4) -- (1,4) -- (1,4.5);
\end{scope}
\node at (10,2.5) {$\mapsto$};
\node[tableau,nodes in empty cells=false] at (11,2)
{
&&&\\
4\\
2&3\\
1&1&2&4\\
};
\end{tikzpicture}
\caption{From a five-vertex lattice path configuration to a SSYT.}
\label{fig:lp2SSYT}
\end{figure}
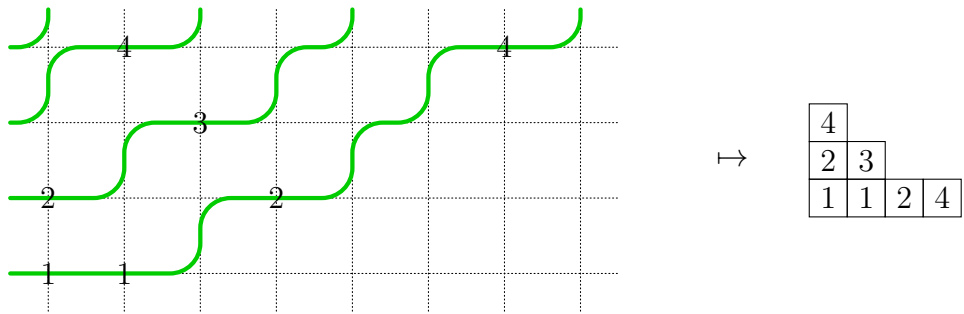

Secondly, five-vertex model configurations with boundary conditions given by $\lambda$
(extended or not -- the extension is unique at $\alpha=0$) are in bijection with SSYTs of shape $\lambda$.
For each path (from bottom to top), record the rows where the path goes East and turn this into one row of the tableau, see Figure~\ref{fig:lp2SSYT}
for an example. Since a weight $z_i$ is assigned to each such East move on row $i$, this clearly matches with formula \eqref{eq:SSYT}.

{\em Remark.} Such five-vertex configurations are also in bijection with {\em pipe dreams}\/ of Grassmannian permutations
\cite{BB-rcgraph,FK-Schubert}, see \cite[\S5.2.5]{hdr}.

\subsection{Cauchy identity}\label{ssec:cauchy}
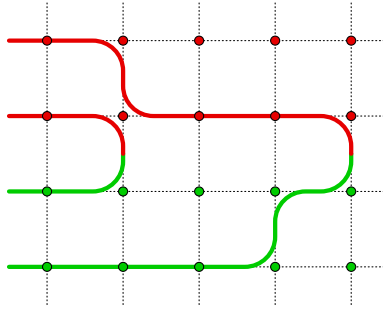
\begin{figure}
\begin{tikzpicture}[line cap=round]
\draw[empty] (0.5,0.5) grid (5.5,4.5);
\begin{scope}[every path/.style={draw=mygreen,ultra thick,rounded corners=4mm}]
\draw (0.5,1) -- (4,1) -- (4,2) -- (5,2) -- (5,2.5);
\draw (0.5,2) -- (2,2) -- (2,2.5);
\end{scope}
\begin{scope}[every path/.style={draw=myred,ultra thick,rounded corners=4mm}]
\draw (0.5,3) -- (2,3) -- (2,2.5);
\draw (0.5,4) -- (2,4) -- (2,3) -- (5,3) -- (5,2.5);
\end{scope}
\foreach\x in {1,2,3,4,5}
{
\foreach\y in {1,2}
\node[part,fill=mygreen] at (\x,\y) {};
\foreach\y in {3,4}
\node[part,fill=myred] at (\x,\y) {};
}
\end{tikzpicture}
\caption{A configuration contributing to $s_\lambda(w_1,\ldots,w_m)s_\lambda(z_1,\ldots,z_{n})$ with $\lambda=(3,1)$.}
\label{fig:cauchy}
\end{figure}

A classical identity is that
\begin{equation}\label{eq:cauchy}
\sum_\lambda s_\lambda(w_1,\ldots,w_m) s_\lambda(z_1,\ldots,z_{n}) = \frac{1}{\prod_{i=1}^m \prod_{j=1}^{n} (1-w_iz_j)}
\end{equation}
(for purposes of convergence, we assume $|w_iz_j|<1$ for all $i,j$).
It is tempting to try to interpret the l.h.s.\ as a partition function. 
An example of a configuration contributing to such a partition function can
be found on Figure~\ref{fig:cauchy}. For simplicity, we use the five-vertex model of the previous section (though the same reasoning would
work with the more general free-fermionic six-vertex model).
The boundary conditions are empty on all sides except
the West side on which $m$ red paths and $n$ green paths enter. We mark vertices with a coloured dot indicating which type of path can go through them:
the green vertices have the usual five-vertex Boltzmann weights,
i.e., a weight of $z_j$ for each step to the right on green row $j$, whereas the configurations
around red vertices are obtained by vertical flip from the usual ones, with a weight of $w_i$
per step to the right on red row $i$. The width $p$ needs to be taken to infinity in order to accommodate for arbitrarily large partitions.

We still have the Yang--Baxter equation for green vertices, and similarly for red vertices, but also a mixed YBE:
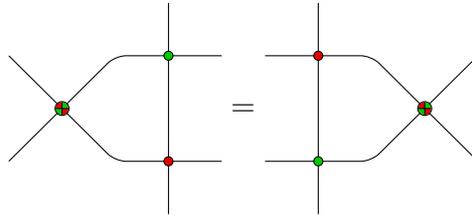
\[
\begin{tikzpicture}[baseline=-3pt,scale=.7]
\draw[rounded corners] (-1,-1) -- (1,1) -- (3,1);
\draw[rounded corners] (-1,1) -- (1,-1) -- (3,-1);
\draw (2,-2) -- (2,2);
\node[part,fill=mygreen] at (2,1) {};
\node[part,fill=myred] at (2,-1) {};
\node[part,circle cross split,circle cross split part fill={myred,mygreen,mygreen,myred}] at (0,0) {};
\end{tikzpicture}
=
\begin{tikzpicture}[baseline=-3pt,scale=.7]
\draw[rounded corners] (-1,-1) -- (1,-1) -- (3,1);
\draw[rounded corners] (-1,1) -- (1,1) -- (3,-1);
\draw (0,-2) -- (0,2);
\node[part,fill=mygreen] at (0,-1) {};
\node[part,fill=myred] at (0,1) {};
\node[part,circle cross split,circle cross split part fill={myred,mygreen,mygreen,myred}] at (2,0) {};
\end{tikzpicture}
\]
where the extra multicoloured crossing is given by the Boltzmann weights
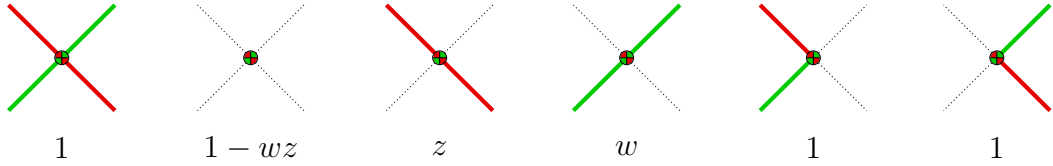
\[
\begin{tikzpicture}
\matrix[column sep=1cm,row sep=5pt,every cell/.style={scale=.7,ultra thick}] {
\draw[mygreen] (-1,-1) -- (1,1);
\draw[myred] (-1,1) -- (1,-1);
\node[part,circle cross split,circle cross split part fill={myred,mygreen,mygreen,myred}] at (0,0) {};
&
\draw[empty] (1,1) -- (-1,-1) (1,-1) -- (-1,1);
\node[part,circle cross split,circle cross split part fill={myred,mygreen,mygreen,myred}] at (0,0) {};
&
\draw[empty] (1,1) -- (-1,-1);
\draw[myred] (1,-1) -- (-1,1);
\node[part,circle cross split,circle cross split part fill={myred,mygreen,mygreen,myred}] at (0,0) {};
&
\draw[mygreen] (1,1) -- (-1,-1);
\draw[empty] (1,-1) -- (-1,1);
\node[part,circle cross split,circle cross split part fill={myred,mygreen,mygreen,myred}] at (0,0) {};
&
\draw[mygreen] (-1,-1) -- (0,0);
\draw[myred] (-1,1) -- (0,0);
\draw[empty] (0,0) -- (1,1) (0,0) -- (1,-1);
\node[part,circle cross split,circle cross split part fill={myred,mygreen,mygreen,myred}] at (0,0) {};
&
\draw[mygreen] (1,1) -- (0,0);
\draw[myred] (1,-1) -- (0,0);
\draw[empty] (0,0) -- (-1,-1) (0,0) -- (-1,1);
\node[part,circle cross split,circle cross split part fill={myred,mygreen,mygreen,myred}] at (0,0) {};
\\
\node{$1$};
&
\node{$1-wz$};
&
\node{$z$};
&
\node{$w$};
&
\node{$1$};
&
\node{$1$};
\\
};
\end{tikzpicture}
\]
(A second mixed YBE can be written with roles of green and red vertices switched, but is trivially related to the first one.)

Let us now show how to obtain \eqref{eq:cauchy} for $m=n=1$.
We have
\begin{alignat*}{2}
&\begin{tikzpicture}[baseline=-3pt]
\draw[mygreen,ultra thick] (-0.5,-.5) -- (0,-.5); \draw[rounded corners] (0,-.5) -- (3.5,-.5) -- (4,0); \draw[empty] (4,0) -- (4.5,-.5);
\draw[myred,ultra thick] (-0.5,.5) -- (0,.5); \draw[rounded corners] (0,.5) -- (3.5,.5) -- (4,0); \draw[empty] (4,0) -- (4.5,.5);
\foreach\x in {0,1,2,3} {
\draw (\x,-.5) -- (\x,.5);
\draw[empty] (\x,-1) -- (\x,-.5) node[part,fill=mygreen] {} (\x,.5) node[part,fill=myred] {} -- (\x,1);
\draw [decorate,decoration={brace}] (0,1.2) -- node[above=1mm] {$p$} (3,1.2);
\node[part,circle cross split,circle cross split part fill={myred,mygreen,mygreen,myred}] at (4,0) {};
}
\end{tikzpicture}
&&=(1-wz)\begin{tikzpicture}[baseline=-3pt]
\draw[mygreen,ultra thick] (-0.5,-.5) -- (0,-.5); \draw (0,-.5) -- (3,-.5); \draw[empty] (3,-.5) -- (3.5,-.5);
\draw[myred,ultra thick] (-0.5,.5) -- (0,.5); \draw (0,.5) -- (3,.5); \draw[empty] (3,.5) -- (3.5,.5);
\foreach\x in {0,1,2,3} {
\draw (\x,-.5) -- (\x,.5);
\draw[empty] (\x,-1) -- (\x,-.5) node[part,fill=mygreen] {} (\x,.5) node[part,fill=myred] {} -- (\x,1);
\draw [decorate,decoration={brace}] (0,1.2) -- node[above=1mm] {$p$} (3,1.2);
}
\end{tikzpicture}
+
\begin{tikzpicture}[baseline=-3pt]
\draw[mygreen,ultra thick] (-0.5,-.5) -- (3.5,-.5);
\draw[myred,ultra thick] (-0.5,.5) -- (3.5,.5);
\foreach\x in {0,1,2,3} {
\draw[empty] (\x,-.5) -- (\x,.5);
\draw[empty] (\x,-1) -- (\x,-.5) node[part,fill=mygreen] {} (\x,.5) node[part,fill=myred] {} -- (\x,1);
\draw [decorate,decoration={brace}] (0,1.2) -- node[above=1mm] {$p$} (3,1.2);
}
\end{tikzpicture}
\\
&&&=(1-wz)\sum_{k=0}^{p-1} s_{(k)}(w) s_{(k)}(z) + (wz)^p
\\
&=\begin{tikzpicture}[baseline=-3pt]
\draw[mygreen,ultra thick] (-1.5,-.5) -- (-1,0); \draw[rounded corners] (-1,0) -- (-.5,-.5) -- (0,-.5) -- (3,-.5);  \draw[empty] (3,-.5) -- (3.5,-.5);
\draw[myred,ultra thick] (-1.5,.5) -- (-1,0); \draw[rounded corners] (-1,0) -- (-.5,.5) -- (0,.5) -- (3,.5); \draw[empty] (3,0.5) -- (3.5,.5);
\foreach\x in {0,1,2,3} {
\draw (\x,-.5) -- (\x,.5);
\draw[empty] (\x,-1) -- (\x,-.5) node[part,fill=myred] {} (\x,.5) node[part,fill=mygreen] {} -- (\x,1);
}
\draw [decorate,decoration={brace}] (0,1.2) -- node[above=1mm] {$p$} (3,1.2);
\node[part,circle cross split,circle cross split part fill={myred,mygreen,mygreen,myred}] at (-1,0) {};
\end{tikzpicture}
&&=\begin{tikzpicture}[baseline=-3pt]
\draw[empty] (-.5,-.5) -- (0,-.5) -- (3,-.5) -- (3.5,-.5);
\draw[empty] (-.5,.5) -- (0,.5) -- (3,0.5) -- (3.5,.5);
\foreach\x in {0,1,2,3} {
\draw[empty] (\x,-.5) -- (\x,.5);
\draw[empty] (\x,-1) -- (\x,-.5) node[part,fill=myred] {} (\x,.5) node[part,fill=mygreen] {} -- (\x,1);
}
\draw [decorate,decoration={brace}] (0,1.2) -- node[above=1mm] {$p$} (3,1.2);
\end{tikzpicture}
=1
\end{alignat*}
where in the last line, we used the fact that green (resp.\ greered) paths can only go North/East (resp.\ South/East).
This allows to compute the sum in the second line, and
by taking $p$ to infinity, one recovers \eqref{eq:cauchy}. The case of general $m$ and $n$ can be treated similarly, by introducing 
a 45 degree rotated $m\times n$ block of multicoloured crossings to the right
of the partition function representing the l.h.s.\ of \eqref{eq:cauchy}, noting that up to vanishingly small terms as $p\to\infty$
one obtains $\prod_{i,j}(1-w_iz_j)$ times the l.h.s., then
applying the YBE repeatedly and finding that the only option is the empty configuration, with weight $1$.

\subsection{Product rule}\label{ssec:puzzle}
Because $\Lambda_n$ is a ring, given a basis such as the $s_\lambda$, one can define structure constants
(which turn out to be independent of $n$):
\[
s_\lambda s_\mu = \sum_{\nu} c_{\lambda,\mu}^\nu s_\nu
\]
These are the celebrated {\em Littlewood--Richardson}\/ coefficients, for which a plethora of 
formulae exists. For example,
$s_{(1)} ^2 =  s_{(2)} + s_{(1,1)}$. In general, they are known to be nonnegative and enumerate various combinatorial
objects such as Littlewood--Richardson tableaux.

One can ask whether methods from integrability can be used to compute the $c_{\lambda,\mu}^\nu$. The answer is
that one can, but on condition that one use a model of rank 2 (i.e., with two conserved quantities at each vertex)
\cite{artic46}, in line with the general philosophy that to multiply and expand two partition functions, the model
should have the sum of ranks of the models associated to the two partition functions \cite{artic71}.

Consider a model on a size $n$ triangular region of a honeycomb lattice (so, with $n^2$ vertices).
The boundary conditions are the ones associated
to $\lambda$ using green paths on the NorthWest side, the ones associated to $\mu$ using red paths on the NorthEast side,
and a new one on the South side, namely that green paths have ending locations given by $\nu$ as usual, except
the remanining spots are taken by red paths.
The allowed vertices are:
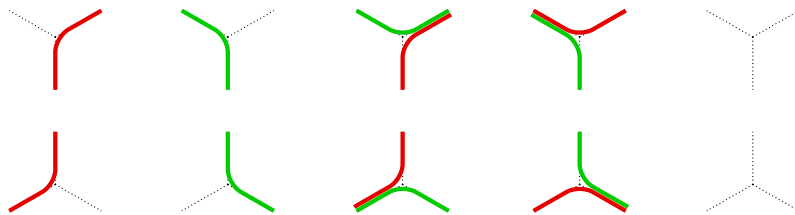
\[
\begin{tikzpicture}
\matrix[every cell/.style={scale=.7,ultra thick,rounded corners=2mm},column sep=1cm,row sep=.5cm,row 2/.style={scale=-1}] {
\draw[empty] (0,0) -- ++(270:1) (0,0) -- ++(150:1) (0,0) -- ++(30:1);
\draw[myred] (270:1) -- (0,0) -- (30:1);
&
\draw[empty] (0,0) -- ++(270:1) (0,0) -- ++(150:1) (0,0) -- ++(30:1);
\draw[mygreen] (270:1) -- (0,0) -- (150:1);
&
\draw[empty] (0,0) -- ++(270:1) (0,0) -- ++(150:1) (0,0) -- ++(30:1);
\draw[myred] (270:1) -- (0,-.1)  -- ++(30:1.05);
\draw[mygreen] (30:1) -- (0,0) -- (150:1);
&
\draw[empty] (0,0) -- ++(270:1) (0,0) -- ++(150:1) (0,0) -- ++(30:1);
\draw[mygreen] (270:1) -- (0,-.1)  -- ++(150:1.05);
\draw[myred] (30:1) -- (0,0) -- (150:1);
&
\draw[empty] (0,0) -- ++(270:1) (0,0) -- ++(150:1) (0,0) -- ++(30:1);
\\
\draw[empty] (0,0) -- ++(270:1) (0,0) -- ++(150:1) (0,0) -- ++(30:1);
\draw[myred] (270:1) -- (0,0) -- (30:1);
&
\draw[empty] (0,0) -- ++(270:1) (0,0) -- ++(150:1) (0,0) -- ++(30:1);
\draw[mygreen] (270:1) -- (0,0) -- (150:1);
&
\draw[empty] (0,0) -- ++(270:1) (0,0) -- ++(150:1) (0,0) -- ++(30:1);
\draw[myred] (270:1) -- (0,-.1)  -- ++(30:1.05);
\draw[mygreen] (30:1) -- (0,0) -- (150:1);
&
\draw[empty] (0,0) -- ++(270:1) (0,0) -- ++(150:1) (0,0) -- ++(30:1);
\draw[mygreen] (270:1) -- (0,-.1)  -- ++(150:1.05);
\draw[myred] (30:1) -- (0,0) -- (150:1);
&
\draw[empty] (0,0) -- ++(270:1) (0,0) -- ++(150:1) (0,0) -- ++(30:1);
\\
};
\end{tikzpicture}
\]
(note that both red and green paths are conserved across vertices).
Then $c_{\lambda,\mu}^\nu$ is the number of configurations with such rules (i.e., all Boltzmann weights are $1$).

A few comments are in order. The model above was actually introduced earlier in the physics
literature in the context of {\em random tilings}:
in \cite{WidM}, Widom showed its equivalence to a model of square-triangle tilings and proved its integrability.
Secondly, there is an easy bijection between the configurations of this model and Knutson--Tao {\em puzzles}
\cite{KT}, whose enumeration is known to reproduce Littlewood--Richardson coefficients.
The bijection consists in replacing every edge with a label $0$, $1$, $10$ according to:
(see also \cite{Pur} for a bijection between puzzles and square-triangle tilings, completing the loop of equivalences)
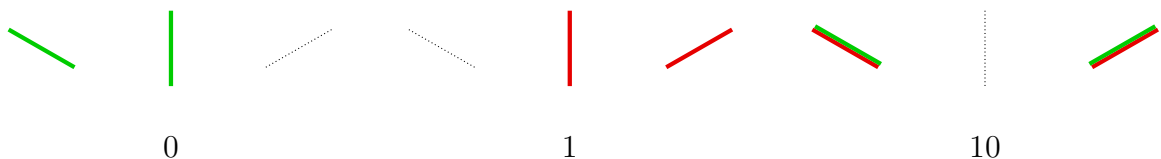
\[
\begin{tikzpicture}
\matrix[column sep=1cm,row sep=.5cm]{
\draw[ultra thick,mygreen] (150:-.5) -- ++(150:1);
&
\draw[ultra thick,mygreen] (270:-.5) -- ++(270:1);
&
\draw[empty] (30:-.5) -- ++(30:1);
&
\draw[empty] (150:-.5) -- ++(150:1);
&
\draw[ultra thick,myred] (270:-.5) -- ++(270:1);
&
\draw[ultra thick,myred] (30:-.5) -- ++(30:1);
&
\draw[ultra thick,myred] (150:-.5) -- ++(150:1);
\draw[ultra thick,mygreen,yshift=1.25pt,xshift=1pt] (150:-.5) -- ++(150:1);
&
\draw[empty] (270:-.5) -- ++(270:1);
&
\draw[ultra thick,myred] (30:-.5) -- ++(30:1);
\draw[ultra thick,mygreen,yshift=1.25pt,xshift=-1pt] (30:-.5) -- ++(30:1);
\\
&\node{$0$};&
&
&\node{$1$};&
&
&\node{$10$};&
\\
};
\end{tikzpicture}
\]

For example, here are the configurations contributing
to $s_{(1)} ^2$ and the corresponding puzzles:
\begin{center}
\begin{tikzpicture}[yscale=-1.732,scale=0.707,rotate=45]
\node at (2.5,-.5) {$(1)$};
\node at (-.5,2.5) {$(1)$};
\node at (2.25,2.25) {$(2)$};
\clip (0,0) -- (0,4) -- (4,0);
\foreach\i in {0,...,3}
\foreach\j in {0,...,3}
{
\coordinate (u-\i-\j) at (\i+1/3,\j+1/3);
\coordinate (d-\i-\j) at (\i+2/3,\j+2/3);
\draw[empty] (u-\i-\j) -- ++(1/3,1/3) (u-\i-\j) -- ++(1/3,-2/3) (u-\i-\j) -- ++(-2/3,1/3);
}
\begin{scope}[ultra thick,mygreen,rounded corners]
\draw (u-0-3) ++(-2/3,1/3) -- (u-0-3) -- (d-0-3);
\draw (u-0-1) ++(-2/3,1/3) -- (u-0-1) -- (d-0-1) -- ([shift={(-.02,-.02)}]u-1-1) -- ++(1/6,-1/3)  ++(.04,.04) -- ([shift={(.02,.02)}]d-1-0) -- (u-2-0) -- ([shift={(-.02,-.02)}]d-2-0) -- ++(1/3,-1/6) ++(.04,.04) -- ([shift={(.02,.02)}]u-3-0) -- (d-3-0);
\end{scope}
\begin{scope}[ultra thick,myred,rounded corners]
\draw (u-3-0) ++(1/3,-2/3) -- ([shift={(-.02,-.02)}]u-3-0) -- ++(-1/3,1/6) ++(.04,.04) -- ([shift={(.02,.02)}]d-2-0) -- (u-2-1) -- (d-2-1);
\draw (u-1-0) ++(1/3,-2/3) -- (u-1-0) -- ([shift={(-.02,-.02)}]d-1-0) -- ++(-1/6,1/3) ++(.04,.04) -- ([shift={(.02,.02)}]u-1-1) -- (d-1-1) -- (u-1-2) -- (d-1-2);
\end{scope}
\end{tikzpicture}\qquad\begin{tikzpicture}[xscale=-1,yscale=-1.732,scale=0.707,rotate=45]
\node at (2.5,-.5) {$(1)$};
\node at (-.5,2.5) {$(1)$};
\node at (2.25,2.25) {$(1,1)$};
\clip (0,0) -- (0,4) -- (4,0);
\foreach\i in {0,...,3}
\foreach\j in {0,...,3}
{
\coordinate (u-\i-\j) at (\i+1/3,\j+1/3);
\coordinate (d-\i-\j) at (\i+2/3,\j+2/3);
\draw[empty] (u-\i-\j) -- ++(1/3,1/3) (u-\i-\j) -- ++(1/3,-2/3) (u-\i-\j) -- ++(-2/3,1/3);
}
\begin{scope}[ultra thick,myred,rounded corners]
\draw (u-0-3) ++(-2/3,1/3) -- (u-0-3) -- (d-0-3);
\draw (u-0-1) ++(-2/3,1/3) -- (u-0-1) -- (d-0-1) -- ([shift={(-.02,-.02)}]u-1-1) -- ++(1/6,-1/3)  ++(.04,.04) -- ([shift={(.02,.02)}]d-1-0) -- (u-2-0) -- ([shift={(-.02,-.02)}]d-2-0) -- ++(1/3,-1/6) ++(.04,.04) -- ([shift={(.02,.02)}]u-3-0) -- (d-3-0);
\end{scope}
\begin{scope}[ultra thick,mygreen,rounded corners]
\draw (u-3-0) ++(1/3,-2/3) -- ([shift={(-.02,-.02)}]u-3-0) -- ++(-1/3,1/6) ++(.04,.04) -- ([shift={(.02,.02)}]d-2-0) -- (u-2-1) -- (d-2-1);
\draw (u-1-0) ++(1/3,-2/3) -- (u-1-0) -- ([shift={(-.02,-.02)}]d-1-0) -- ++(-1/6,1/3) ++(.04,.04) -- ([shift={(.02,.02)}]u-1-1) -- (d-1-1) -- (u-1-2) -- (d-1-2);
\end{scope}
\end{tikzpicture}
\\
\begin{tikzpicture}[execute at begin picture={\bgroup\tikzset{every path/.style={}}\clip (-2.48,-.43201) rectangle ++(4.96,4.46413);\egroup},x={(.93630cm,0cm)},y={(0cm,-.93630cm)},baseline=(current  bounding  box.center),every path/.style={draw=black,fill=none},line join=round]
null\begin{scope}[every path/.append style={fill=white,line width=.01873cm}]
\path[fill=mygray] (-.5,.86603) -- (0,0) -- (.5,.86603) -- cycle;
\path[fill=mygray] (-.5,.86603) -- (.5,.86603) -- (0,1.73205) -- cycle;
\node[black] at (.25,.43301) {0};
\node[black] at (-.25,.43301) {1};
\node[black,,scale=.83565] at (0,.86603) {10};
\path[fill=mygray] (0,1.73205) -- (.5,.86603) -- (1,1.73205) -- cycle;
\path[fill=mygray] (0,1.73205) -- (1,1.73205) -- (.5,2.59808) -- cycle;
\node[black] at (.75,1.29904) {1};
\node[black] at (.25,1.29904) {1};
\node[black] at (.5,1.73205) {1};
\path[fill=mygray] (.5,2.59808) -- (1,1.73205) -- (1.5,2.59808) -- cycle;
\path[fill=mygray] (.5,2.59808) -- (1.5,2.59808) -- (1,3.46410) -- cycle;
\node[black] at (1.25,2.16506) {0};
\node[black] at (.75,2.16506) {0};
\node[black] at (1,2.59808) {0};
\path[fill=mygray] (1,3.46410) -- (1.5,2.59808) -- (2,3.46410) -- cycle;
\node[black] at (1.75,3.03109) {1};
\node[black,,scale=.83565] at (1.25,3.03109) {10};
\node[black] at (1.5,3.46410) {0};
\path[fill=mygray] (-1,1.73205) -- (-.5,.86603) -- (0,1.73205) -- cycle;
\path[fill=mygray] (-1,1.73205) -- (0,1.73205) -- (-.5,2.59808) -- cycle;
\node[black] at (-.25,1.29904) {0};
\node[black] at (-.75,1.29904) {0};
\node[black] at (-.5,1.73205) {0};
\path[fill=mygray] (-.5,2.59808) -- (0,1.73205) -- (.5,2.59808) -- cycle;
\path[fill=mygray] (-.5,2.59808) -- (.5,2.59808) -- (0,3.46410) -- cycle;
\node[black,,scale=.83565] at (.25,2.16506) {10};
\node[black] at (-.25,2.16506) {0};
\node[black] at (0,2.59808) {1};
\path[fill=mygray] (0,3.46410) -- (.5,2.59808) -- (1,3.46410) -- cycle;
\node[black] at (.75,3.03109) {1};
\node[black] at (.25,3.03109) {1};
\node[black] at (.5,3.46410) {1};
\path[fill=mygray] (-1.5,2.59808) -- (-1,1.73205) -- (-.5,2.59808) -- cycle;
\path[fill=mygray] (-1.5,2.59808) -- (-.5,2.59808) -- (-1,3.46410) -- cycle;
\node[black] at (-.75,2.16506) {0};
\node[black] at (-1.25,2.16506) {1};
\node[black,,scale=.83565] at (-1,2.59808) {10};
\path[fill=mygray] (-1,3.46410) -- (-.5,2.59808) -- (0,3.46410) -- cycle;
\node[black] at (-.25,3.03109) {1};
\node[black] at (-.75,3.03109) {1};
\node[black] at (-.5,3.46410) {1};
\path[fill=mygray] (-2,3.46410) -- (-1.5,2.59808) -- (-1,3.46410) -- cycle;
\node[black] at (-1.25,3.03109) {0};
\node[black] at (-1.75,3.03109) {0};
\node[black] at (-1.5,3.46410) {0};
\end{scope}
\end{tikzpicture}\quad\begin{tikzpicture}[execute at begin picture={\bgroup\tikzset{every path/.style={}}\clip (-2.48,-.43201) rectangle ++(4.96,4.46413);\egroup},x={(.93630cm,0cm)},y={(0cm,-.93630cm)},baseline=(current  bounding  box.center),every path/.style={draw=black,fill=none},line join=round]
null\begin{scope}[every path/.append style={fill=white,line width=.01873cm}]
\path[fill=mygray] (-.5,.86603) -- (0,0) -- (.5,.86603) -- cycle;
\path[fill=mygray] (-.5,.86603) -- (.5,.86603) -- (0,1.73205) -- cycle;
\node[black] at (.25,.43301) {0};
\node[black] at (-.25,.43301) {1};
\node[black,,scale=.83565] at (0,.86603) {10};
\path[fill=mygray] (0,1.73205) -- (.5,.86603) -- (1,1.73205) -- cycle;
\path[fill=mygray] (0,1.73205) -- (1,1.73205) -- (.5,2.59808) -- cycle;
\node[black] at (.75,1.29904) {1};
\node[black] at (.25,1.29904) {1};
\node[black] at (.5,1.73205) {1};
\path[fill=mygray] (.5,2.59808) -- (1,1.73205) -- (1.5,2.59808) -- cycle;
\path[fill=mygray] (.5,2.59808) -- (1.5,2.59808) -- (1,3.46410) -- cycle;
\node[black] at (1.25,2.16506) {0};
\node[black] at (.75,2.16506) {1};
\node[black,,scale=.83565] at (1,2.59808) {10};
\path[fill=mygray] (1,3.46410) -- (1.5,2.59808) -- (2,3.46410) -- cycle;
\node[black] at (1.75,3.03109) {1};
\node[black] at (1.25,3.03109) {1};
\node[black] at (1.5,3.46410) {1};
\path[fill=mygray] (-1,1.73205) -- (-.5,.86603) -- (0,1.73205) -- cycle;
\path[fill=mygray] (-1,1.73205) -- (0,1.73205) -- (-.5,2.59808) -- cycle;
\node[black] at (-.25,1.29904) {0};
\node[black] at (-.75,1.29904) {0};
\node[black] at (-.5,1.73205) {0};
\path[fill=mygray] (-.5,2.59808) -- (0,1.73205) -- (.5,2.59808) -- cycle;
\path[fill=mygray] (-.5,2.59808) -- (.5,2.59808) -- (0,3.46410) -- cycle;
\node[black] at (.25,2.16506) {1};
\node[black,,scale=.83565] at (-.25,2.16506) {10};
\node[black] at (0,2.59808) {0};
\path[fill=mygray] (0,3.46410) -- (.5,2.59808) -- (1,3.46410) -- cycle;
\node[black] at (.75,3.03109) {0};
\node[black] at (.25,3.03109) {0};
\node[black] at (.5,3.46410) {0};
\path[fill=mygray] (-1.5,2.59808) -- (-1,1.73205) -- (-.5,2.59808) -- cycle;
\path[fill=mygray] (-1.5,2.59808) -- (-.5,2.59808) -- (-1,3.46410) -- cycle;
\node[black] at (-.75,2.16506) {1};
\node[black] at (-1.25,2.16506) {1};
\node[black] at (-1,2.59808) {1};
\path[fill=mygray] (-1,3.46410) -- (-.5,2.59808) -- (0,3.46410) -- cycle;
\node[black] at (-.25,3.03109) {0};
\node[black] at (-.75,3.03109) {0};
\node[black] at (-.5,3.46410) {0};
\path[fill=mygray] (-2,3.46410) -- (-1.5,2.59808) -- (-1,3.46410) -- cycle;
\node[black,,scale=.83565] at (-1.25,3.03109) {10};
\node[black] at (-1.75,3.03109) {0};
\node[black] at (-1.5,3.46410) {1};
\end{scope}
\end{tikzpicture}
\end{center}

\gdef\MRshorten#1 #2MRend{#1}%
\gdef\MRfirsttwo#1#2{\if#1M%
MR\else MR#1#2\fi}
\def\MRfix#1{\MRshorten\MRfirsttwo#1 MRend}
\renewcommand\MR[1]{\relax\ifhmode\unskip\spacefactor3000 \space\fi
\MRhref{\MRfix{#1}}{{\scriptsize \MRfix{#1}}}}
\renewcommand{\MRhref}[2]{%
\href{http://www.ams.org/mathscinet-getitem?mr=#1}{#2}}
\bibliographystyle{amsalphahyper}
\bibliography{biblio}
\end{document}